# A COURSE OF ALGEBRA (PART I)

**Frasser C.E.**
*Ph.D. In Engineering Sciences, Odessa National Polytechnic University, Ukraine*

**Abstract**
The following is an exposition of a course of algebra that Prof. Aleksandr Aleksandrovich Zykov (1922-2013) distributed among the participants of his seminar in graph theory not far away from Odessa, Ukraine, on September, 1991. It is a privilege for me to be able to reproduce, with some good additions, the English version of this remarkable course that he designed carefully over several years with the help of other colleagues and with the only purpose of making the science of algebra more accessible to the less gifted student. Prof. Zykov was an admirable person and a very clever mathematician and scientist. His scientific interests were wide: from the mathematical sciences at any level and in any of their branches to the theory of relativity as I recall. He kindly wrote a review recommending the results of my doctoral dissertation, which was necessary for the approval to obtain a Ph.D. degree during the meeting of my defense. In that meeting, I was very nervous because I had in front of me fifteen professors all of them doctors of technical or mathematical sciences asking all sort of questions and whom I had to convince of the fact that my results were true. I am sure this algebra course played the role that Prof. Zykov, with his deepness of thought, had primarily intended, that is, one course designed with the right approach to help out many of his students.

## A.A. ZYKOV
## LECTURES ON ALGEBRA

ABSTRACT
A systematic, but not entirely traditional algebra course for university's mathematical faculties, designed for students with far from ideal preparation in high school, mainly taught by the author for four years at the Samarkand University named after Alisher Navoi. The introductory lectures were written jointly with Associate Professor Asad Amonovich Amonov and Senior Lecturer Yusuf Sharipovich Sharipov (Department of Algebra and Number Theory).

### PREFACE

One of my friends decided to become a trumpeter, but he never got into music school. The admissions committee found out he had no ear nor spirit. However, let's imagine he still ended up as a student at the school (due to a shortage or simply "on a call"): what would musicians do with him later? after all, often, in such an absurd position, teachers of the mathematical faculty of the deputy state university, where I was lucky enough to work between 1982 and 1987, found themselves. I said I was lucky without irony. Although it was eventually possible to get rid of students who had absolutely no relation to mathematics, I still had to deal, for the most part, with those whose hearing was not sufficiently developed and whose spirit was not yet set at all. Thus, if it were not for the warm, friendly attitude and support of colleagues, then, under such difficult conditions, this course of algebra could not have appeared despite the fact that it was designed for students with only the most necessary minimum of mathematical abilities and very incomplete school knowledge. But it can be recommended as a starting point even for the students of Moscow State University. The presented material in it needs only additions, but not revision.

We rejected both the great amount of focusing required for the widely used textbook of A.G. Kurosh and the excessive academicism of the textbook written by A.I. Kostrikin. At the same time, I recommend it for independent study of some sections and for additional reading. The presentation of the material is fully consistent. The task of students' assimilation of the very spirit of algebra as the science of algebraic operations is absolutely relevant and the difference between the algebraic approach to the study of mathematical objects and the analytical one is explained already in the lectures and it is emphasized throughout the course.

I strongly disagreed with my colleagues at the department when, due to the weak composition of the students, they avoided the issue of algebraic expansion in their lectures. After all, it means to teach algebra without algebra, which is natural and logical in technical colleges where you have to solve mathematical problems by all available means and not to cultivate separately "the spirit of algebra", "the spirit of geometry" and "the spirit of analysis". This is absolutely unacceptable for future mathematicians both as researchers and teachers. If from the field of rational numbers, the field of real numbers arises analytically, then its next extension to the field of complex numbers arises algebraically. A fact you can't hide, how in vain one hero tried in O. Henry: "He denied the allegations, but he could not deny the alligators" (play on words from the novel "The man higher up": Figure it out!).

I am strongly against the exclusion from the compulsory curriculum of such a topic as solving the equations of third and fourth degree: you can't turn a course of algebra into "Ivan who doesn't understand any kind of family connection or affinity". At the same time, I consider it is entirely justified to limit ourselves to the formulation of the root existence theorem, which, on its own, is the fundamental theorem of algebra. The proof of the latter is purely analytical in nature and it is unfeasible within the framework of algebra itself.

A systematic study of the equations with one unknown begins with the cases of first and second degree not only to eliminate school shortcomings, but also to provide full justification: after all, everything has to be done in the complex field. By the way, attempts to fill the gaps with a refresher course in school mathematics usually turned out to be ineffective. It is unlikely that in 1-2 months you will be able to master well all of that for which you did not do much in the course of 10 years. Here, it is more likely that a completely different approach will help out. At first, I should captivate students with unusual material not only based on ready-made knowledge, but also made with certain minimum of intelligence and common sense. Only then, in the process of going through something new and interesting, they gradually do something "right" to make up for doing something wrong.

There is no place for the theory of determinants in the first semester: for algebra itself, there they are not needed yet, and for analytic geometry, the information on third-order determinants that directly follow from the properties of the mixed product of three vectors is quite sufficient. But in the second semester, $n$th order determinants have the full right to appear not just as computing tools, but as multilinear functions of the elements of a linear space, which also simplifies the presentation. In the same way, why to study in the first semester (before the general theory of linear spaces) the Gauss-Jordan elimination method, and especially start the course with it. In some cases when it is necessary to solve a system of linear equations, it is possible to get by with school methods (a reminder of this is given in the historical plan of the introductory lecture and can be reinforced in the very first practical lesson). Finally, symmetric polynomials (such a beauty should not be thrown out of the course completely) are considered in the second semester in connection with the general concept of algebras.

At first glance, the theory of operations on matrices at the end of the first semester may seem out of place. However, in fact, determinants are not needed to decide degeneracy (solution exists and is not unique) and non-degeneracy (solution exists and is unique) and to find the inverse matrix. Rings of a square matrix naturally fit in the general study of algebraic systems while the addition and multiplication of rectangular matrices, on the one hand, allow for a number of clear and important interpretations, and on the other hand, they are also studied in the general terms of algebraic operations. Such a transfer of different topics allows to avoid in the second semester overloads due to the theory of determinants and symmetric polynomials.

The material of the second semester, that is, linear algebra and multidimensional geometry with the exception of the previously mentioned modifications, is traditional. The third semester consists of two parts "homomorphisms and factorization" and "permutation groups" including elements of combinatorial analysis and the famous Polya theorem. In our opinion, knowledge of the latter is as necessary as that of number theory, for which a separate course is set aside in the second semester.

The book is written in outline form and is not a manual for self-education. If a weak student dares to read it alone, then the most one can hope for is an illusion of understanding the lecture. However, the one who overcomes the basics under the guidance of an experienced teacher will be able to consciously work through the material of the order of the topic, let's say, "complex numbers" (where, in particular, it becomes inevitable to turn to other teaching aids). As the experience shows, for a weak, but diligent student, the decisive criteria for whether he can continue his studies at the Faculty of Mathematics is its ability to fully understand the theorem on the uniqueness of the neutral element of a groupoid.

In addition to the material on the geometric representation and properties of complex numbers, requests for additional sources require tasks for independent work (three such tasks in the first semester are provided) an also familiarity with the theorem on partitioning a set into equivalence classes; the latter, in our opinion, should be presented in the course "Introduction to speciality area program".

AUXILIARY REFERENCES for the algebra course in the first semester:
A.G. Kurosh, Course of Higher Algebra.
A.I. Kostrikin, Introduction to Algebra.
PROBLEM BOOKS:
D.K. Faddeev, I.S. Sominskiy, Collection of Problems in Higher Algebra.
I.B. Proskuryakov, Collection of Problems in Linear Algebra.

A.A. Amonov, A.A. Zykov, Yu. Sh. Sharipov

FROM THE ART OF SOLVING EQUATIONS TO THE SCIENCE OF ALGEBRAIC OPERATIONS
(Introductory Lecture on Algebra)

As is evident from the surviving written documents of Ancient Egypt and Babylonia, already 3,000 years BC, not only four arithmetic operations on integers and fractional numbers (as well as extracting the root) were known, but also the solution of some problems of the type when based on the result of a series of operations on an unknown number, it was necessary to establish this number. The simplest of those problems leads to a first-degree equation, which in modern notation has one of the forms

$$x \pm a = b, \quad a \pm x = b, \quad ax = b.$$

(In fact, the problem and its solution were written in words and specific numbers since there were no letter symbols, no mathematical operation signs and no equal sign. But for the sake of convenience in conveying the essence of the matter, we will agree to use the now familiar symbols). To solve the third equation, Egyptians immediately divided $b$ by $a$ (in a more complicated way than we do, though) while Babylonians first found the reciprocal value for $a$ using a pre-compiled table and then multiplied it by $b$.

More complex problems of the type

$$ax + bx = c$$

were not solved according to the modern scheme $(a+b)x = c$, $x = c/(a+b)$ because equivalent transformations of equations, in particular, the reduction of like terms containing the desired number were not then known. In such cases, they used the *false position method*, which we will explain using one problem from the Egyptian papyrus: the quantity and its fourth part together give 15. "Count this for me", that is, find the original quantity.

The solution began with the following words: "count from 4; from it, you must take a quarter, namely 1; together 5". Then division 15/5 = 3 and multiplication 4(3) = 12 were performed, giving the desired number. In other words, to solve the equation $x + \frac{1}{4}x = 15$, first we set $x = 4$ (in a simpler way so that no fractions appear), as a result of which the left side will take the value of 5, and since this is three times less than 15, to obtain the "true $x$" from the "false one" the latter must be tripled. This path seems awkward to us, but for a purely arithmetic solution of problems leading to systems of two first-degree equations with two unknowns, the false position method still works (see exercise 1).

In Babylonia, more complex problems were solved, for example, those leading to systems of the type

$$x + y = a,$$
$$xy = b.$$

The formulation of the problem itself looked something like this: "The length and width together are equal to $a$, the area is $b$; what are the length and width?" (It is clear that we are talking about a rectangle). The solution process was written down as a recitation: "break off half of $a$, you will get $a/2$" (all, of course, using specific numbers). "Square $a/2$, subtract $b$ and take the root of the difference" (there were tables for the last step as well as for finding reciprocals). "By adding the resulting number to $a/2$, you will find the length and by subtracting it from $a/2$, - the width". Finally, there was a check that the numbers found satisfied the conditions of the problem, in other words, that "you solved it correctly".

The solution procedure corresponds to the following formulas:

$$x = \frac{a}{2} + \sqrt{\left(\frac{a}{2}\right)^2 - b}, \qquad y = \frac{a}{2} - \sqrt{\left(\frac{a}{2}\right)^2 - b},$$

but who personally and, most importantly, how it was discovered is unknown: the prescription-dogmatic presentation ("do it this way") does not contain any evidence or even explanations of why it should be done this way and not otherwise. Checking at the end (Egyptians also did it for some simple problems) can be considered only the embryo of a proof.

Of particular interest is the way in which the Chinese around 2,000 BC dealt with problems of the type

$$a_1 x + b_1 y + c_1 z = d_1$$
$$a_2 x + b_2 y + c_2 z = d_2$$
$$a_3 x + b_3 y + c_3 z = d_3,$$

there was no notation for the desired quantities, nor addition sign, multiplication sign, equal sign, and the system was written in the form of a table of specific numbers

| $a_1$ | $a_2$ | $a_3$ |
|---|---|---|
| $b_1$ | $b_2$ | $b_3$ |
| $c_1$ | $c_2$ | $c_3$ |
| $d_1$ | $d_2$ | $d_3$, |

over the columns of which actions were then performed corresponding to the sequential elimination of unknowns. All numbers were positive (others were not known at that time), however, in the process of solving the problem, it was often necessary to subtract a larger number from a smaller one. Then, they did the opposite, they subtracted the smaller from the larger, but wrote the difference in a different color of ink. Thus, the same numbers had to be operated depending on their color: the black number had to be subtracted where the red one had to be added, and vice versa. However, further steps that prepared the introduction of negative numbers took place (somewhere else) much later, and even the highly developed mathematics of ancient Greece passed this by.

One of the greatest achievements of the Greek scientists of the 6th – 2th centuries BC was the development of logic and its wide application to the natural sciences, which made possible, in particular, the transformation of mathematics from a set of isolated facts and rules into a systematic science. Euclid's famous "Elements" (around 300 BC) opened up the formulation of the basic principles referred as axioms and postulates, from which then, with the help of logical reasoning, the theorems of arithmetic and geometry were derived. Algebra appeared in geometric form; known and unknown quantities were not numbers, but segments, rectangles, parallelepipeds (in fact, their lengths, areas, volumes). With this approach, for example, equation $ax = b$ should be written as

$$ax = b^2$$

(so that the dimensions of both sides coincide), and the formulation of the problem looked like this: "on a given segment, construct a rectangle equal in area to a given square". More complex problems on the application of a rectangle to a segment "with a deficiency" and "with an excess" (all three problems and their connection with conic sections; parabola, ellipse and hyperbola, should be discussed in the course of analytic geometry and more detailed information can be found in textbooks on the history of mathematics) led to equations of the form

$$ax \pm x^2 = b^2,$$

which were solved by geometric constructions with full justification. The conditions for the existence of a solution to this (or that) equation were called *diorisms* - they responded with a modern inequality between the coefficients that ensured the non-negativity of the discriminant and the positivity of at least one of the roots (see exercise 5).

The reason that the main characters of ancient Greek algebra were not numbers, but geometric objects, was the discovery in the 4th century BC by the Pythagoreans of such an unexpected phenomenon as incommensurability, in particular, the fact that the length of the diagonal of a unit square cannot be expressed by either a whole or a fractional number (if, of course, we mean an exact value, but not an approximate one). Let's give a *proof by contradiction* (again in modern notation) following the Pythagorean tradition.

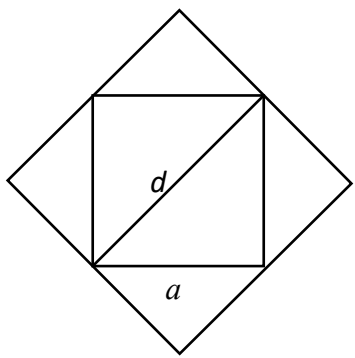

Fig. 1

Let's assume that the diagonal of a square with side $a = 1$ has length $d$ (Fig. 1). Since the area of the square constructed on this diagonal having side length $d$ is twice the area of the original square $a^2 = 1$, $d^2 = 2$. Number $d$ cannot be a positive integer because $1^2 < 2$, but already $2^2 > 2$, furthermore $3^2 > 2$, and so on. Let $d$ be a fractional number. We will show that this assumption leads to a contradiction.

Let's denote by $m/n$ an irreducible fraction, that is, a fraction in simplest terms, which is equal to $d$. Then $(m/n)^2 = 2$ meaning that

$$m^2 = 2n^2.$$

This implies that $m$ is even, otherwise $m^2$ as the product of two identical odd numbers, would itself be odd, whereas $2n^2$ is even. Writing $m$ in the form $2k$, we obtain from $m^2 = 2n^2$ the expression $4k^2 = 2n^2$ or after dividing by two

$$n^2 = 2k^2,$$

as above, it follows that $n$ is even.

And so, the numbers $m$ and $n$ are both even and the fraction $m/n$ can be reduced by dividing both of its terms by two, which contradicts the assumption of its irreducibility, that is, the assumption that when writing $d$ as $m/n$, all possible reductions to simplest terms have already been made.

Excessive attention to logical perfection at the expense of the applied side of mathematics, associated with direct measurements and calculations (it was considered an occupation worthy of slaves and artisans, but not of the learned elite, and the greatest mathematician, physicist, and engineer of antiquity, Archimedes, was one of the few exceptions) had its negative consequences, which later served as one of the main reasons for the decline of ancient Greek science. In relation to arithmetic and algebra, this was manifested in the fact that irrational numbers were never discovered by the Greeks (there was a very subtle and profound theory for the relationships of incommensurable quantities, but these relationships themselves were not considered numbers) and diorisms created a strong immunity against the emergence of negative and imaginary numbers.

Mathematics of India and the Arab countries of the first millennium AD, having absorbed the achievements of the ancient science of Egypt, Babylonia, China (partly) and Greece, gave the world an important innovation: the decimal positional system and the number zero. Although the positional principle itself was widely used in ancient Babylonia, the base of this system (60 instead of 10) made it much more cumbersome to calculate (the remains of the Babylonian system have been preserved in the form of dividing the hour into minutes and seconds, and the circle into degrees, arc minutes and seconds). Indian scholars of the 6th-12th centuries also mention negative numbers, with an explanation of arithmetic operations on them, but only as a "logical possibility", denying these numbers equality of rights with the positive ones.

Although many problems of algebraic content were solved before our era and also during the eight centuries of our era, algebra as a separate science did not yet exist, and its founder is rightfully considered to be the outstanding Central Asian scientist al-Khwarizmi (783-850). He and his followers developed general methods that allow, without changing the roots of an equation, to bring it to a form in which these roots are obvious or at least easier to find than in the original one.

The main types of transformation in al-Khwarizmi are *Al-Jabr* ("restoration") and *Al-Muqabala* ("comparison"). The first consists of removing the subtracted term from one side of the equation with subsequent restoration as a positive term in the other side; the second means the reduction of like terms. To these transformations is added the previously known division (as well as multiplication) of both sides by the same number. For example, to solve linear equation $7x - 27 = 9 - 5x$, a chain of transformations is performed

$$7x - 27 = 9 - 5x,$$
$$5x + 7x - 27 = 9,$$
$$5x + 7x = 9 + 27,$$
$$12x = 36,$$
$$x = 3.$$

Each line, including the last one, is an equation *equivalent* to the original one in the sense that any specific number $x$ that satisfies one equation also satisfies the other. Each of the questions of existence, uniqueness and actual calculation of the root has the same answer for all equations. But if it is not possible to "guess" the root of the original equation right away (and the question of uniqueness will not be resolved by this), then for the last equation all indicated problems are trivial: its only root, according to the pertinent remark of the algebraists of those times, is known to us whether one wishes to or not.

Using transformations, al-Khwarizmi reduces any equation of first or second degree to one of the six canonical (simplest) forms, where each degree of the unknown can be contained no more than once and only as addend, not as subtrahend (see exercise 7) and for each case he gives a geometric solution. From "Al-Jabr" Europeans subsequently formed the word *algebra*, and from the name of "al-Khwarizmi" the term *algorithm*, meaning a general method (a single system of rules) for solving any specific problem from a certain class of problems.

To the merits of al-Khwarizmi as a mathematician, it should be added that he was one of the first to appreciate the decimal positional number system and actively promoted it. In his arithmetical treatise, he clearly sets out the rules for operating on this system, including column multiplication and the use of zero. This subsequently led to the emergence and widespread use of decimal fractions in the East earlier than in Europe (in the observatory of Ulugbek, Uzbekistan's great astronomer, sines and tangents of angles were calculated with an accuracy of up to eight decimal places).

Using the algebraic techniques of al-Khwarizmi, outstanding mathematician (as well as philosopher and poet) of the 11th-12th centuries Omar Khayyam reduces any equation of degree 1, 2 or 3 to one of the 25 canonical forms (see exercise 7 again), which then solves using geometric algebra and the ancient Greek theory of conic sections. But its solution is closer in spirit to finding the intersection point of the graphs of two functions rather than constructing a segment. Only an unfortunate omission in the analysis of one of the cases prevented Khayyam from coming to the conclusion that a cubic equation can have three different roots (see exercises 8 and 9). In general, scientists from the countries of the Near and Middle East treated the works of their predecessors with deep respect; in particular, they translated and carefully studied the works of ancient Greek philosophers and mathematicians, and in the presentation of known results, and in obtaining new ones, the ancient rigor of proofs was preserved, but at the same time serious attention was paid to the application of mathematical discoveries in the real world (road and canal construction, architecture, land surveying, trade, complex problems of property division, and so on), as well as to other sciences (optics, geography and especially astronomy). Irrational numbers received the "right of citizenship" since they can always be replaced with decimal fractions with the necessary accuracy and attempts were made (fundamentally interesting, but mathematically, not completed, unfortunately) to build a unified theory of rational and irrational numbers.

When speaking about the emergence of algebra as a special branch of mathematics and the successes of the algebraists of the medieval East, it would be wrong to pass over in silence some omissions. We have already mentioned one of them, related to cubic equations. Another fact, more general, concerns negative numbers: instead of bringing this "logical possibility" of the Indian scholars to systematic development and application, they further aggravated the tradition of avoiding such quantities (and information about the "red" and "black" numbers of the Chinese either did not reach Central Asia at all, or was not properly assessed there either). "Intolerance" of negative numbers explains the abundance of canonical forms of equations in al-Khwarizmi and Omar Khayyam. It was also, apparently, one of the main reasons why Khayyam did not create (long before Descartes and Fermat) analytic geometry, even though, many historians believe that he stood on the threshold of this discovery.

In Europe, the lush flowering of ancient science and culture gave way to a long medieval hibernation. The mathematics of this period made almost no major discoveries; it was mainly a process of accumulation of material which, while important for the future, was very slow. Only from the 14th century did interest in the works of mathematicians of Greece, India and Central Asia awaken, in particular, both mathematical treatises of al-Khwarizmi were translated from Arabic into Latin. One of the merits of the mathematicians of the late European Middle Ages should be considered a more tolerant attitude towards negative quantities than in Asia: in addition to the Indian interpretation ("debt" as opposed to "property"), they were given the now well-known kinematic explanation. Despite this, even later, in the Renaissance, it was still a long way off from the equality of the importance of both signs, and algebraists, working with negative numbers, called them "fictitious" (an end to this discrimination was put only in the 19th century).

The new rise of science in the Renaissance period, of course, also embraced algebra. Europeans of the 16th century were especially interested in solving equations. In Italy, there were even public competitions to find exact numerical values (and not the construction of corresponding segments) for the roots of an equation of degree higher than 2. Del Ferro, Tartaglia, Cardano and Ferrari discovered general methods for the algebraic solution of equations of third and fourth degrees. It was found that in a number of cases it was impossible to avoid intermediate actions with expressions of the form $\sqrt{-1}$, to which at the same time no meaning could be attributed and, therefore, were called "sophistical" or "imaginary" numbers.

As noted by Cardano, equation

$$x^2 - 10x + 40 = 0,$$

according to ancient Greek diorism, has no solution, but it is satisfied by the values $x_1 = 5 + \sqrt{-15}$ and $x_2 = 5 - \sqrt{-15}$ if we act on the roots of negative numbers according to the same rules as on the roots of positive numbers, in particular, we consider $(\sqrt{-15})^2 = -15$ and the founder of the literal algebraic calculus, the French mathematician Viète, discovered that even for such numbers as $x_1$ and $x_2$, his theorem on the connection between the coefficients and roots of an equation remains valid.

Only at the turn of the 18th and 19th centuries did expressions of the form $a + b\sqrt{-1}$ with reals $a$ and $b$ get the name *complex numbers*, and their representation by points of the plane eliminated the mysticism (the decisive role in this process, which then led to the complete recognition of complex numbers, belongs to Gauss). But before that, it was necessary to act with such numbers without understanding their meaning (nature) and guided only by the rules established for "understandable" numbers. The situation was similar with regard to "dumb" (irrational) numbers, which could not be expressed exactly by any fraction, although the possibility of their simple geometric representation on a straight line and approximate representation by a decimal fraction with any required accuracy weakened the mystical touch.

The geometric algebra of the ancient Greeks inevitably relied on the letter notation of known and unknown quantities and, in general, the replacement of specific numbers with letters for setting and describing the course of solving a problem was used by Jordanes Nemorarius (12th-13th centuries). But his approach did not yet mean the creation of a letter computation since the result of the action on two expressions was designated each time by a new letter. For example, the law $(a + b) c = ac + bc$ looked like this: "If *d* is the sum of *a* and *b*, *e* is the product of *d* and *c*, *f* is the product of *a* and *c*, *g* is the product of *b* and *c*, *h* is the sum of *f* and *g*, then *e* is equal to *h*."

Viète's symbolism was not distinguished by its compactness either: he wrote our equality $A^3 + 3BA = D$ in the form

*A* cubus + *B* planum in *A*3 aequatur D solido;

however, the origins of algebraic calculus are located precisely here because such a notation allows one to act with the literal expressions themselves as with numbers. The notation was soon improved by Harriot, whose notation of the same equality

$$aaa + 3ba = d$$

is much closer to the modern one. Around the same time (the second half of the 16th century), along with the signs +, −, × (the latter was often omitted), the = sign appeared, and the first letter of the Latin word *radix* (root) in combination with a line above $\sqrt{\phantom{x}}$, which then replaced brackets, formed the radical sign in Viète.

The emergence of letter algebra greatly strengthened the determination to act with various kinds of quantities, including those "incomprehensible" according to certain rules.

Here we will interrupt the historical narrative (further information of this nature is better given during the algebra course itself and in the course of history of mathematics) and will only note a characteristic feature of modern algebra. If the *analytical* approach to the study of certain objects is associated primarily with the clarification of their internal structure, then with the *algebraic* approach, the laws by which these objects are combined with each other come to the fore. Mathematical examples will be sufficient later, but for now we will illustrate the importance of the algebraic approach (and the dialectical unity of both approaches) using an example that is not related to algebra itself.

When discovering the periodic law, D.I. Mendeleev did not proceed from the internal structure of atoms (then still unknown) but from the "external" properties of chemical elements, and primarily, their ability to form certain compounds with each other. However, the periodic table itself played a decisive role in the subsequent creation of the atomic model to find the mass and charge of nuclei and the number of electrons in the shells of a particular chemical element. Thus, the original "algebraic" approach opened the doors wide to the analytical approach.

EXERCISES

1. Read carefully the story by A.P. Chekhov "The Tutor" and solve the problem of buying blue and black cloth the way little Pete should have done, that is, starting with a "false position": let's assume that everything bought was blue…

2. Write down the solution of the system

$$x - y = a,$$
$$xy = b.$$

using modern formulas and in the form of a Babylonian recipe.

*Solution.* (Try it yourself first).  $x = \sqrt{\left(\frac{a}{2}\right)^2 + b} + \frac{a}{2},$   $y = \sqrt{\left(\frac{a}{2}\right)^2 + b} - \frac{a}{2}$

3. In ancient Babylonia, more complex systems were also solved, for example

$$x \pm y = a,$$
$$x^2 + y^2 = b.$$

by reducing to those considered above. Write down the part of the recipe that reduces the solution of the system (separately for the cases of + sign and - sign in the first equation) to the system given in the lecture text, and "docking" this part with the already known recipe.

*Solution.* Case 1:  $x = \frac{a}{2} + \sqrt{\frac{b}{2} - \left(\frac{a}{2}\right)^2},$   $y = \frac{a}{2} - \sqrt{\frac{b}{2} - \left(\frac{a}{2}\right)^2}$

Case 2:  $x = \sqrt{\frac{b}{2} - \left(\frac{a}{2}\right)^2} + \frac{a}{2},$   $y = \sqrt{\frac{b}{2} - \left(\frac{a}{2}\right)^2} - \frac{a}{2}$

4. Write down the system of equations

$$x + 2y + 2z = 9,$$
$$2x + 5y + z = 17,$$
$$2x + 7y + 2z = 22$$

in the form of a Chinese table, on which you carry out a successive elimination of unknowns in order to find their values. Before the negative term put also a + sign, but provide the term itself with a line on top (meaning "another color").

*EXPLANATION.* To eliminate $x$ from the second and third equations, we now subtract twice the first from these equations. In Chinese notation, the corresponding operations on the columns of the table of coefficients (and constant terms) transform it as follows:

$$\begin{array}{ccc} 1 & 2 & 2 \\ 2 & 5 & 7 \\ 2 & 1 & 2 \\ \hline 9 & 17 & 22 \end{array} \longrightarrow \begin{array}{ccc} 1 & 0 & 0 \\ 2 & 1 & 3 \\ 2 & 3 & \bar{2} \\ \hline 9 & \bar{1} & 4 \end{array}$$

then, using similar actions, but leaving the second column unchanged, we will give the second row the form  0   1   0. Finally, using the third column (having previously divided all its numbers by 7), we will give the third row the form  0   0   1. The numbers under the line, which arose as a result of all these operations with the columns, will be the sought values of the unknowns $x, y, z$.

4'. It is more natural for us to write the coefficients (and free terms) of the equations in the same order as the equations themselves, that is, in rows, not in columns. The system under consideration then corresponds to the table

$$\begin{array}{ccc|c} 1 & 2 & 2 & 9 \\ 2 & 5 & 1 & 17 \\ 2 & 7 & 2 & 22 \end{array}$$

and the solution process in general comes down to the proper application of three types of operations on the rows:

a) Replacing a row with its sum (or difference) with another row multiplied by an arbitrary number: the "other row" itself does not change;
b) multiplication or division of all numbers in a row by the same non-zero number;
c) Interchanging rows (not used in our example).

Solve the system of exercise 4 using the entries just described. Instead of a line above the number, you can now put a minus in front, as is the modern practice.

4". Try to apply the solution method described in exercise 4' to the systems

$$\begin{aligned} x + 2y + z &= 4, \\ 2x + y + 2z &= 5, \\ 3x + 3y + 3z &= 8, \end{aligned} \qquad \begin{aligned} x + y + z &= 3, \\ x + 2y + z &= 4, \\ 2x + y + 2z &= 5 \end{aligned}$$

and explain the results.

5. From the point of view of geometric algebra, how many solutions does each of the following equations have: $2x - x^2 = 1$, $3x - x^2 = 3$, $x^2 + 3x + 2 = 0$, $ax + x^2 = b$? Could all these equations have been found among the ancient Greeks?

6. Prove the irrationality of numbers $\sqrt{3}$ and $\sqrt{2} + \sqrt{3}$.

6'. Can the sum of irrational numbers be a rational number?

7. Explain the number of canonical equations in al-Khwarizmi (6) and in Omar Khayyam (25), taking into account that equations, which obviously did not have positive roots, were not considered.

8. (For those who are especially interested). Read the book by A.P. Yushkevich "History of Mathematics in the Middle Ages" (Moscow, 1961) about O. Khayyam's approach to solving cubic equations and find out what the "unfortunate omission" is.

9. The equation

$$(x-1)(x-2)(x-3) = 0$$

is obviously of third degree and has three different roots. Why do you think simple examples of this kind did not lead Omar Khayyam to the idea of the existence of three roots in a cubic equation?

10. Solve the equation $x^2 + x + 1 = 0$ and check the validity of Viète's theorem for the obtained expressions (complex numbers).

## PART 1. BASIC ALGEBRAIC STRUCTURES: GROUPS, RINGS AND FIELDS

## MAPPINGS

Let $X$ and $Y$ be two non-empty sets. Their *Cartesian product* $X \times Y$ is the set

$$\{(x,y)/x \in X, y \in Y\},$$

that is, the set of all ordered pairs of elements, the first of which belongs to set $X$ and the second to set $Y$.

*Example 1.* If $X = \{a, b, c\}, Y = \{1, 2\}$, then
$X \times Y = \{(a,1), (a,2), (b,1), (b,2), (c,1), (c,2)\}$,
$Y \times X = \{(1,a), (1,b), (1,c), (2,a), (2,b), (2,c)\} \neq X \times Y$,
and if $a, b, c \notin \{1, 2\}$, then $(X \times Y) \cap (Y \times X) = \emptyset$.

*Example 2.* If $X = \{a, b\}$, $Y = \{b, c\}$, where $a \neq c$, then
$X \times Y = \{(a,b), (a,c), (b,b), (b,c)\}$,
$Y \times X = \{(b,a), (b,b), (c,a), (c,b)\} \neq X \times Y$,
but $(X \times Y) \cap (Y \times X) = \{(b,b)\} \neq \emptyset$.

*Example 3.* If $X$ and $Y$ are mutually perpendicular numerical axes and the point $O$ of their intersection serves as the origin on each of them, then $X \times Y$ is the Cartesian coordinate plane $xOy$.

When $Y = X$, the product $X \times X$ is denoted by $X^2$ and is called the *Cartesian square* of set $X$. For example, if $X = \{a, b\}$, then $X^2 = \{(a,a), (a,b), (b,a), (b,b)\}$ and all four pairs are different. Note that ordered pairs $(p, q)$ and $(r, s)$ are considered to be the same, that is, $(p, q) = (r, s)$ if and only if both their first and second components coincide: $p = r$ and $q = s$. (Unordered pairs of ordered pairs $\{p, q\}$ and $\{r, s\}$ are the same, if $p = r$, $q = s$ or $p = s$, $q = r$. Pair $\{p, q\}$, $p \neq q$, is a two-element set, but when $p = q$, the set is the singleton one $\{p\}$). Similarly, the Cartesian cube $X^3 = \{(x, y, z)/x, y, z \in X\}$ is the set of all ordered triplets of elements $X$, and so on. It is natural to consider by definition that $X^1 = X$ and $X^0 = \emptyset$.

*Question.* Is it true that $X^2 \times X = X \times X^2 = X^3$ and $X \times X^0 = X$?

Let $M \neq \emptyset$ and $N \neq \emptyset$ be given. The *mapping* $f: M \to N$ of a set $M$ onto a set $N$ assigns to each element $x \in M$ a certain element $f(x) \in N$. $f$ is also called a *function* defined on a set $M$ and taking values in a set $N$. This mapping can be strictly defined as a subset of the set of pairs $M \times N$, in which for each $x \in M$, there exists a unique pair $(x, y)$ and then its second element $y = f(x) \in N$ is called the *image* of the element $x$ under the mapping $f$ (or the value of the function $f$ for the value $x$ of its argument).

A mapping $f: M \to N$ is called *injective* (or one-to-one) if different elements of $M$ have different images in $N$, that is,

$$\forall x_1, x_2 \in M: x_1 \neq x_2 \Rightarrow f(x_1) \neq f(x_2).$$

A mapping $f: M \to N$ is called *surjective* if each element $y \in N$ has at least one preimage (meaning such an element $x \in M$, whose image is $y$), that is, if

$$\forall y \in N \, \exists x \in M: f(x) = y.$$

A mapping $f: M \to N$ is *bijective* (is a *bijection*, or a one-to-one and at the same time surjective correspondence between $M$ and $N$) if it is simultaneously injective and surjective.

*Example 4.* In Fig. 2, $f_1(a) = 3$, $f_1(b) = 3$, $f_1(c) = 4$, $f_1(d) = 4$, $f_1(e) = 3$.

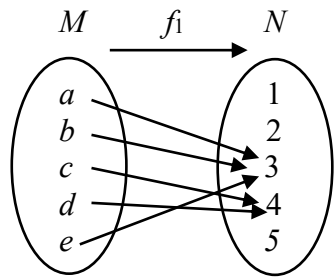

Fig. 2

The mapping $f_1 : M \to N$ is neither injective nor surjective.

*Example 5.* In Fig. 3, $f_2(a) = 1$, $f_2(b) = 2$, $f_2(c) = 4$, $f_2(d) = 5$, $f_2(e) = 6$.

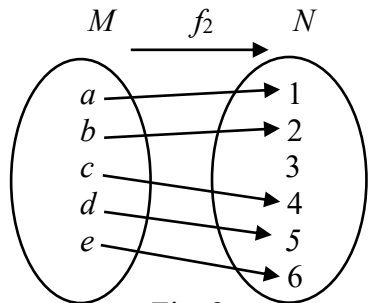

Fig. 3

The mapping $f_2 : M \to N$ is injective but not surjective: for element $3 \in N$ there is no $x \in M$ such that $f_2(x) = 3$.

*Example 6.* In Fig. 4, $f_3(Vanya) = Lida$, $f_3(Petya) = Lida$.

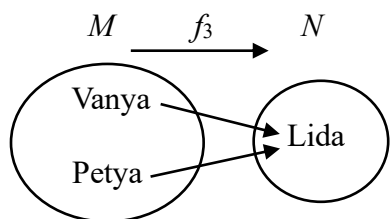

Fig. 4

The mapping $f_3 : M \to N$ is surjective but not injective.

*Example 7.* In Fig. 5, the mapping $f_4 : M \to N$ is injective and surjective, that is, it is a bijection of the set $M$ onto the set $N$.

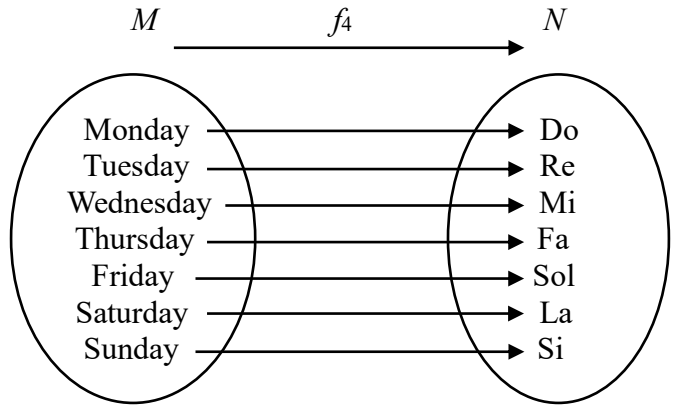

Fig. 5

*Example 8.* But this is not a mapping at all (Fig. 6) (elements $ and * from *M* do not correspond to anything in *N*); it is only a correspondence relating some elements from *M* to some from *N*.

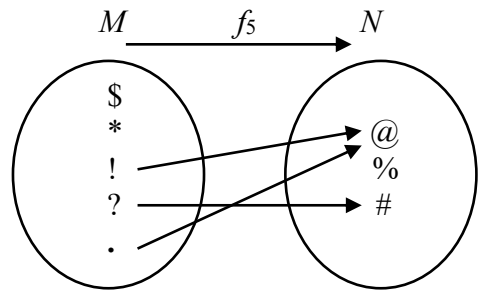

Fig. 6

*Problem.* For each of the functions $f_1(x) = |x|$, $f_2(x) = x^2$, $f_3(x) = x^3$, $f_4(x) = \sqrt{x}$, $f_5(x) = -\sqrt{x}$, $f_6(x) = \sqrt{|x|}$, specify its own pair of sets $M_i, N_i \subseteq \mathbb{R}$ (that is, as subsets of real numbers) so that the mapping $f_i : M_i \to N_i$ is a) injective, b) surjective, c) bijective ($i = 1, 2,\ldots,6$).

Let $f(x, y)$ be a function of two arguments (e.g. $x + y$, $x - y$, $x \cdot y$ for $x, y \in \mathbb{R}$). $f$ can be viewed as a mapping (a function of one argument) that assigns to each ordered pair $(x, y)$ some element:
$$f : M^2 \to N.$$

In general, a function $f(x_1, x_2, \ldots, x_m)$ of *m* arguments, which relates to each ordered set of *m* elements of a set *M* some element of a set *N*, is a mapping
$$f : M^m \to N \quad (m = 0, 1, 2,\ldots).$$

Number *m* is called the *arity* of mapping *f*.

0-ary (nullary) mapping fixes some element in *N*, usually a constant (since here the "function" *f* has no arguments at all on which its value could depend);

1-ary (unary) assigns to each $x \in M$ an element of $f(x) \in N$;

2-ary (binary) assigns to each ordered pair of elements $(x, y) \in M$ an element $f(x, y) \in N$;

3-ary (ternary) assigns to each ordered triple $(x, y, z) \in M$ an element $f(x, y, z) \in N$ and so on. When $M = N$, mappings are also called *operations* (nullary, unary, binary, ternary, etc.). In algebra, a binary operation is often denoted in general form by the symbol $*$ and the result of its application to elements $x, y \in M$ (in the given order) is written in the form $x * y$; with this kind of notation (in particular, when replacing the symbol $*$ with the sign of a particular operation, for example, $+, -, \cdot, :$) the operation $*$ itself is usually called a *composition*, and the result of its application to the pair $x, y \in M$, e.g. the element $x * y \in N = M$, is called a composition of elements $x$ and $y$.

A non-empty set with a set of operations defined on it (arbitrariness of arities) is called an *algebraic system*.

## ALGEBRAIC SYSTEMS WITH ONE BINARY OPERATION

A *groupoid* $(M, *)$ is a non-empty set $M$ with a defined composition $*$ on it, which uniquely assigns some element $x * y \in M$ to each ordered pair of elements $x, y \in M$.

*Examples 9.* If $M = \mathbb{N} = \{1, 2, 3, ...\}$ (the set of all natural numbers) and $*$ is the ordinary addition $(+)$, then $(\mathbb{N}, +)$ is a groupoid. $(\mathbb{N}, *)$ is also a groupoid (where $*$ is the ordinary multiplication). At the same time, $(\mathbb{N}, -)$ and $(\mathbb{N}, :)$ are not groupoids (why?). $(\{1, 2, ..., 10\}, +)$ and $(\{1, 2, ..., 10\}, \cdot)$ are not groupoids since the addition and multiplication of the set of the first ten natural numbers can lead beyond its limits: $7 + 8 \notin \{1, 2, ..., 10\}$, $3 \cdot 4 \notin \{1, 2, ..., 10\}$ even though $7, 8, 3, 4 \in \{1, 2, ..., 10\}$. If $x * y$ means $x^y$, then $(\mathbb{N}, *)$ is a groupoid; in it, unlike the groupoids of the previous examples, not always $x * y = y * x$ (For instance, $2^3 \neq 3^2$).

Despite the extreme generality of the definition of a groupoid (and, as a consequence, the almost complete absence of information about groupoids in general), it is possible to introduce a very important definition and prove an equally important theorem of a general nature.

An element $e$ of a groupoid $(M, *)$ is called *neutral* if $x * e = e * x = x$ for any $x \in M$, or in another notation

$$\forall x \in M: x * e = e * x = x.$$

*Remark 1.* Not every groupoid has a neutral element. For instance, in $(\mathbb{N}, +)$, there is no such an element. However, if $M = \mathbb{N}_0 = \{0, 1, 2, 3, ...\}$ (the set of all whole numbers), then in $(\mathbb{N}_0, +)$, there is one: number 0. In $(\mathbb{N}, \cdot)$, there is also one: number 1. In groupoid $(\mathbb{N}, *)$, where $x * y = = x^y$, we have $x * 1 = x$ for any $x$. This implies that number 1 is a *right neutral element*, but not a *left neutral* one because $1 * x \neq x$ (if $x \neq 1$). Therefore, 1 is not a neutral element (and no other natural number is either).

**Theorem 1.** *In a groupoid, there cannot be more than one neutral element.*

*Proof.* Let $e_1$ and $e_2$ be neutral elements of groupoid $(M, *)$. We will show that $e_1 = e_2$.
If $x * e_1 = x$ for any $x \in M$, then, in particular, $e_2 * e_1 = e_2$. If $e_2 * x = x$ also for any $x \in M$, then, in particular, $e_2 * e_1 = e_1$. But composition $e_2 * e_1$ is a uniquely determined element of set $M$. Therefore, $e_2 = e_1$. □

*Remark 2.* In fact, a stronger statement is proven here: *If a groupoid $(M, *)$ has both a right neutral element $e_1$ such that $\forall x \in M: x * e_1 = x$ and a left neutral element $e_2$ such that $\forall x \in M: e_2 * x = x$, then $e_1 = e_2$ and this element (and only this one) serves as a neutral element in $(M, *)$.*

A groupoid $(M, *)$ is called a *semigroup* if the composition $*$ is associative (satisfies the associativity principle or associative property), that is,

$$\forall x, y, z \in M: x * (y * z) = (x * y) * z.$$

*Examples 10.* $(\mathbb{N}_0, +)$ and $(\mathbb{N}, \cdot)$ are semigroups since always

$$x + (y + z) = (x + y) + z \text{ and } x(yz) = (xy)z.$$

*Example 11.* On the contrary, $(\mathbb{N}, *)$ with $x * y = x^y$ is not a semigroup because equality $x^{(y^z)} = (x^y)^z$ is not always true: $3^{(1^2)} \neq (3^1)^2$.

*Example 12.* Groupoid $(\mathbb{Q}, *)$, where $\mathbb{Q}$ is the set of all rational numbers and $x * y = \frac{x+y}{2}$ is the arithmetic mean, is not a semigroup since for any $x, y, z \in \mathbb{Q}$ the values

$$x * (y * z) = \frac{x + (y+z)/2}{2} \text{ and } (x * y) * z = \frac{(x+y)/2 + z}{2}$$

do not coincide.

A semigroup having a neutral element is called a *monoid*. If this neutral element $e$ is explicitly stated, then monoid $(M, *, e)$ is an algebraic system with two operations: binary (composition $*$) and nullary (fixed element $e$); but it is usually written as $(M, *)$ without specifically noting the nullary operation.

For a monoid $(M, *)$, the following definition makes sense: an element $y \in M$ is called the *inverse* of an element $x \in M$ if $x * y = y * x = e$.

**Theorem 2.** *No element of a monoid can have more than one inverse.*

*Proof.* Let an element $x$ of a monoid $(M, *)$ have inverses $y_1$ and $y_2$. We will show that $y_1 = y_2$.

Due to the associativity of the composition $*$, the expression $y_1 * x * y_2$ (without brackets) makes sense since for specific $y_1, x, y_2 \in M$, the resulting element will be the same for both arrangements of brackets:

$$y_1 * x * y_2 \begin{cases} = y_1 * (x * y_2) = y_1 * e = y_1 \\ = (y_1 * x) * y_2 = e * y_2 = y_2 \end{cases} \Rightarrow y_1 = y_2. \qquad \square$$

*Remark 3.* As in the theorem on the uniqueness of the neutral element of a groupoid, here a stronger statement is actually proven: *if an element $x$ of a monoid has a right inverse $y_2$ such that* $x * y_2 = e$ and a left inverse $y_1$ such that $y_1 * x = e$, then $y_1 = y_2$ and this element (and only this one) serves as the inverse of $x$.

*Remark 4.* Not necessarily each element of a monoid is *invertible*, that is, has an inverse. But if $x \in M$ is an invertible element, then its only inverse is denoted by $x^{-1}$. Since by the definition of the inverse element, $x * x^{-1} = x^{-1} * x = e$, the inverse of $x^{-1}$ is $x$, that is, $(x^{-1})^{-1} = x$.

Note that the associativity of composition $*$ implies the meaningfulness of expressions of the form $a * b * c * d * f$ (with any number of "factors") since any arrangement of brackets leads to the same result. The proof of this, in general, is cumbersome, and initially, instead, we prefer to offer the derivation of individual equalities as exercises, for example, prove that

$$((a * b) * c) * (d * f) = a * \left(b * ((c * d) * f)\right)$$

by using the associative property successively.

**Theorem 3.** *The composition of any number of invertible elements of a monoid is invertible.*

*Proof.* Let $a, b, c, d$ be invertible elements of monoid $(M, *)$. We will show that for the composition $a * b * c * d$, the inverse element will be $d^{-1} * c^{-1} * b^{-1} * a^{-1}$ (the same reasoning applies to any number of elements instead of four).

$$(a * b * c * d) * (d^{-1} * c^{-1} * b^{-1} * a^{-1}) = (a * b * c) * (d * d^{-1}) * (c^{-1} * b^{-1} * a^{-1})$$

$$= (a * b * c) * e * (c^{-1} * b^{-1} * a^{-1}) = (a * b * c) * (c^{-1} * b^{-1} * a^{-1}) = \cdots = a * a^{-1} = e,$$

and similarly,

$$(d^{-1} * c^{-1} * b^{-1} * a^{-1}) * (a * b * c * d) = \cdots = (d^{-1} * c^{-1} * b^{-1}) * (b * c * d)$$

$$= \cdots = d^{-1} * d = e. \qquad \square$$

*Remark 5.* And so, *the inverse element for the composition of invertible elements is the composition of inverse elements arranged in reverse order*.

*Analogy.* When you get dressed, you put on an undershirt first, then a shirt, a jacket, a coat. On the other hand, when you undress, you start with the coat, then you take off the jacket, shirt, undershirt (for the sake of decency, we only consider the upper half of the body).

In every monoid, there is one invertible element: this is the neutral element $e$. A monoid, in which all elements are invertible, is called a *group*. In view of the special importance of this algebraic system, we will formulate its definition without explicitly using the words groupoid, semigroup and monoid.

A system $(M, *)$ consisting of a non-empty set $M$ and a binary operation (composition $*$) defined on it is called a *group*, or, in other words, *a set $M$ forms a group with respect to the operation $*$ if the following four conditions are met:*

1. for any $x, y \in M$, the element $x * y \in M$ is uniquely defined;
2. operation $*$ is associative;
3. $M$ has a neutral element;
4. all elements of $M$ are invertible.

*Examples 13.* $(\mathbb{N}, +)$ is not a group since the third condition is not met.

$(\mathbb{Z}, +)$, where $\mathbb{Z}$ is the set of all integers, is a group ($e = 0$, the role of $x^{-1}$ is played by $-x$).

$(\mathbb{Q}, +)$ is a group.

$(\mathbb{Q}, \cdot)$ is a monoid, but not a group since $0^{-1}$ does not exist (the black sheep started to spoil everything).

$(\mathbb{Q} \setminus \{0\}, \cdot)$ is a group. The set of all rational numbers distinct from zero forms a group under multiplication.

**Theorem 4.** *In a group* $(M, *)$ *for any* $a, b \in M$ *each of the equations*

$$a * x = b, \qquad y * a = b$$

*is uniquely solvable.*

*Proof.* Firstly, let's assume that an element $x \in M$, satisfying the first equation, does exist, and let's find out what that element $x$ is (thereby proving its uniqueness). Then, we will check that the element found really satisfies the equation (that is, we will prove its existence). For the second equation, the reasoning is similar and it is supposed to be carried out independently.

So, let $x$ be a specific element of the group for which $a * x = b$. Then

$$a^{-1} * (a * x) = a^{-1} * b \Rightarrow (a^{-1} * a) * x = a^{-1} * b \Rightarrow e * x = a^{-1} * b \Rightarrow x = a^{-1} * b$$

and that is the only way it can be. Let's check that this $x$ satisfies the equation:

$$a * (a^{-1} * b) = (a * a^{-1}) * b = e * b = b. \qquad \square$$

The operation $*$ in any groupoid $(M, *)$ is called *commutative* and the groupoid itself is called commutative if $\forall x, y \in M: x * y = y * x$.

A commutative group is also called an *abelian* group. In all previous examples, the groups were abelian. Examples of non-abelian groups include, in particular, the following ones:

PERMUTATION GROUPS. *A permutation of nth degree* is called a bijection $\propto: \{1, 2, ..., n\} \to \{1, 2, ..., n\}$. If

$$\propto(1) = i_1, \propto(2) = i_2, ..., \propto(n) = i_n$$

($i_1, i_2, \ldots, i_n \in \{1, 2, \ldots, n\}$ and all are different), then permutation $\alpha$ is written in the form

$$\alpha = \begin{pmatrix} 1 & 2 & \ldots & n \\ i_1 & i_2 & \ldots & i_n \end{pmatrix}.$$

The columns of this array can be rearranged: what matters is what each of the numbers $1, 2, \ldots, n$ turns into, not where it is stated.

*Example 14.* If

$$\alpha \begin{matrix} 1 \searrow 1 \\ 2 \nearrow 2 \\ 3 \longrightarrow 3 \end{matrix}, \text{ to } \alpha = \left\{ \begin{matrix} \begin{pmatrix} 1 & 2 & 3 \\ 2 & 1 & 3 \end{pmatrix} = \begin{pmatrix} 1 & 3 & 2 \\ 2 & 3 & 1 \end{pmatrix} = \begin{pmatrix} 2 & 1 & 3 \\ 1 & 2 & 3 \end{pmatrix} \\ = \begin{pmatrix} 2 & 3 & 1 \\ 1 & 3 & 2 \end{pmatrix} = \begin{pmatrix} 3 & 1 & 2 \\ 3 & 2 & 1 \end{pmatrix} = \begin{pmatrix} 3 & 2 & 1 \\ 3 & 1 & 2 \end{pmatrix} \end{matrix} \right\}.$$

*Remark 6.* Note that $\alpha$ turns 1 into 2, 2 into 1, and 3 into 3. The first row of the six permutation arrays previously described determines the six different permutations of 3 elements taken 3 at a time, that is, the six different ways in which 3 people can stand in line, whose number is $3! = 6$, while the second row of each permutation, turns 1 into 2, 2 into 1 and 3 into 3 each time with respect to the first row. Note also that the 6 permutations are all equal.

Let $S_n$ be the set of all $n!$ permutations of $n$th degree,

$$\alpha = \begin{pmatrix} 1 & 2 & \ldots & n \\ i_1 & i_2 & \ldots & i_n \end{pmatrix} \in S_n, \quad \beta = \begin{pmatrix} i_1 & i_2 & \ldots & i_n \\ j_1 & j_2 & \ldots & j_n \end{pmatrix} \in S_n.$$

Then *the product of $\alpha$ and $\beta$* is called the permutation

$$\alpha\beta = \begin{pmatrix} 1 & 2 & \ldots & n \\ j_1 & j_2 & \ldots & j_n \end{pmatrix} \in S_n.$$

*Example 15.* If $n = 3$,

$$\alpha = \begin{pmatrix} 1 & 2 & 3 \\ 1 & 3 & 2 \end{pmatrix}, \quad \beta = \begin{pmatrix} 1 & 2 & 3 \\ 2 & 1 & 3 \end{pmatrix} = \begin{pmatrix} 1 & 3 & 2 \\ 2 & 3 & 1 \end{pmatrix},$$

then

$$\alpha\beta = \begin{pmatrix} 1 & 2 & 3 \\ 1 & 3 & 2 \end{pmatrix} \cdot \begin{pmatrix} 1 & 3 & 2 \\ 2 & 3 & 1 \end{pmatrix} = \begin{pmatrix} 1 & 2 & 3 \\ 2 & 3 & 1 \end{pmatrix}$$

and

$$\beta\alpha = \begin{pmatrix} 1 & 3 & 2 \\ 2 & 3 & 1 \end{pmatrix} \cdot \begin{pmatrix} 1 & 2 & 3 \\ 1 & 3 & 2 \end{pmatrix} = \begin{pmatrix} 1 & 3 & 2 \\ 2 & 3 & 1 \end{pmatrix} \cdot \begin{pmatrix} 2 & 3 & 1 \\ 3 & 2 & 1 \end{pmatrix} = \begin{pmatrix} 1 & 3 & 2 \\ 3 & 2 & 1 \end{pmatrix} \neq \alpha\beta,$$

so, the multiplication of permutations is in general non-commutative.

*Remark 7.* According to the definition of the product of $\alpha$ and $\beta$, the second row of $\alpha$ must coincide with the first row of $\beta$ and the product $\alpha\beta$ is the array whose first row is the first row of $\alpha$ and second row is the second row of $\beta$. Note that now $\alpha$ turns 1 into 1, 2 into 3 and 3 into 2 while $\beta$ turns 1 into 2, 2 into 1 and 3 into 3. Furthermore, according to Example 14, the 2 permutations used for $\beta$ are equivalent as are the second permutations in the product of $\beta\alpha$. Note that when multiplying $\alpha$ by $\beta$, it is not necessary to rearrange the columns in the permutation array $\beta$, you can reason on the original array: "$\alpha$ converts 1 to 1, $\beta$ converts 1 to 2 meaning $\alpha\beta$ converts 1 to 2. On the other hand, $\alpha$ converts 2 to 3, $\beta$ converts 3 to 3 meaning $\alpha\beta$ converts 2 to 3. Finally, $\alpha$ converts 3 to 2, $\beta$ converts 2 to 1 meaning $\alpha\beta$ converts 3 to 1".

**Theorem 5.** *For any $n \in \mathbb{N}$, a set $S_n$ forms a group under the multiplication of permutations. For $n = 1, 2$, these groups are abelian. For $n \geq 3$, they are non-abelian.*

*Proof.* Let's check that all four conditions of the group definition for system $(S_n, \cdot)$ are satisfied.

1. If $\alpha, \beta \in S_n$, then $\alpha\beta \in S_n$, too. This immediately follows from the definition of the product of permutations and, at the same time, is a special case of a more general statement: *let $\alpha: M \to M$ and $\beta: M \to M$ be bijections. Then the mapping $\alpha\beta: M \to M$, which assigns to each $x \in M$ an element $\beta(\alpha(x)) \in M$, is also a bijection.* The last statement can be easily proved by establishing separately the injection and surjection of mapping $\alpha\beta$.

2. Let $\alpha, \beta, \gamma \in S_n$, $i \in \{1, 2, \ldots, n\}$ be any such that $\alpha(i) = j, \beta(j) = k, \gamma(k) = 1$. Then $\beta\gamma(j) = 1$, $\alpha(\beta\gamma)(i) = 1$, $\alpha\beta(i) = k$, $(\alpha\beta)\gamma(i) = 1$. So, no matter what number $i \in \{1, 2, \ldots, n\}$ we take, both permutations $\alpha(\beta\gamma)$ and $(\alpha\beta)\gamma$ will turn it into the same number 1, and this means that $\alpha(\beta\gamma)$ and $(\alpha\beta)\gamma$ are the same permutation from $S_n$, that is, $\alpha(\beta\gamma) = (\alpha\beta)\gamma$.

3. The neutral element in $S_n$ is the *identical permutation* $e = \begin{pmatrix} 1 & 2 & \ldots & n \\ 1 & 2 & \ldots & n \end{pmatrix}$. It is easy to see that for any $\alpha = \begin{pmatrix} 1 & 2 & \ldots & n \\ i_1 & i_2 & \ldots & i_n \end{pmatrix}$, the identity $\alpha e = e\alpha = \alpha$ is true.

4. If $\alpha = \begin{pmatrix} 1 & 2 & \ldots & n \\ i_1 & i_2 & \ldots & i_n \end{pmatrix}$, then $\alpha^{-1} = \begin{pmatrix} i_1 & i_2 & \ldots & i_n \\ 1 & 2 & \ldots & n \end{pmatrix}$. Indeed, $\begin{pmatrix} 1 & 2 & \ldots & n \\ i_1 & i_2 & \ldots & i_n \end{pmatrix} \cdot \begin{pmatrix} i_1 & i_2 & \ldots & i_n \\ 1 & 2 & \ldots & n \end{pmatrix} =$

$= \begin{pmatrix} 1 & 2 & \ldots & n \\ 1 & 2 & \ldots & n \end{pmatrix} = e$ and $\begin{pmatrix} i_1 & i_2 & \ldots & i_n \\ 1 & 2 & \ldots & n \end{pmatrix} \cdot \begin{pmatrix} 1 & 2 & \ldots & n \\ i_1 & i_2 & \ldots & i_n \end{pmatrix} = \begin{pmatrix} i_1 & i_2 & \ldots & i_n \\ i_1 & i_2 & \ldots & i_n \end{pmatrix} = e$.

Group $S_1$ consists of a single element $e = \begin{pmatrix} 1 \\ 1 \end{pmatrix}$, and group $S_2$ consists of two: $e = \begin{pmatrix} 1 & 2 \\ 1 & 2 \end{pmatrix}$ and $\alpha = \begin{pmatrix} 1 & 2 \\ 2 & 1 \end{pmatrix}$, where $ee = e, e\alpha = \alpha e = \alpha, \alpha\alpha = e$, and both groups are obviously abelian. The two permutations $\alpha, \beta \in S_3$ of example 15, such that $\alpha\beta \neq \beta\alpha$, can easily be adapted to prove the non-commutativity of any $S_n$ with $n \geq 3$.

Let $\alpha = \begin{pmatrix} 1 & 2 & 3 & 4 \dots n \\ 1 & 3 & 2 & 4 \dots n \end{pmatrix}, \beta = \begin{pmatrix} 1 & 2 & 3 & 4 \dots n \\ 2 & 1 & 3 & 4 \dots n \end{pmatrix} \in S_n$, then

$\alpha \beta = \begin{pmatrix} 1 & 2 & 3 & 4 \dots n \\ 1 & 3 & 2 & 4 \dots n \end{pmatrix} \cdot \begin{pmatrix} 1 & 3 & 2 & 4 \dots n \\ 2 & 3 & 1 & 4 \dots n \end{pmatrix} = \begin{pmatrix} 1 & 2 & 3 & 4 \dots n \\ 2 & 3 & 1 & 4 \dots n \end{pmatrix}$,

$\beta \alpha = \begin{pmatrix} 1 & 2 & 3 & 4 \dots n \\ 2 & 1 & 3 & 4 \dots n \end{pmatrix} \cdot \begin{pmatrix} 2 & 1 & 3 & 4 \dots n \\ 3 & 1 & 2 & 4 \dots n \end{pmatrix} = \begin{pmatrix} 1 & 2 & 3 & 4 \dots n \\ 3 & 1 & 2 & 4 \dots n \end{pmatrix} \neq \alpha \beta$. □

MULTIPLICATIVE AND ADDITIVE NOTATION OF A GROUP. For an arbitrary group $(M, *)$, a *multiplicative* notation is often used: composition $*$ is denoted by $\cdot$ (and this dot is sometimes omitted), the neutral element - by 1 (the inverse of $x$ as before - by $x^{-1}$) although set $M$ is not necessarily numerical and composition may differ from ordinary multiplication. An abelian group is often written *additively*, denoting $*$ by $+$, the neutral element by 0, and the inverse of $x$ by $-x$, calling it the *opposite*.

ABSTRACT GROUPS. In all the specific examples of groups considered above, the approach to their study was mainly analytical: based on the nature of the elements, the law of their composition was established. With a purely algebraic approach, the nature of the elements is not important, only the law of composition itself is important. It can be specified, for instance, by a *Cayley table*, which, for the group $(M, *)$ with $M = \{x_1, x_2, \dots, \}$, looks like this (table 1):

Table 1

| * | $x_1$ | $x_2$ | ... | $x_j$ | ... |
|---|---|---|---|---|---|
| $x_1$ | $x_1 * x_1$ | $x_1 * x_2$ | ... | $x_1 * x_j$ | ... |
| $x_2$ | $x_2 * x_1$ | $x_2 * x_2$ | ... | $x_2 * x_j$ | ... |
| ... | ... | ... | ... | ... | ... |
| $x_i$ | $x_i * x_1$ | $x_i * x_2$ | ... | $x_i * x_j$ | ... |
| ... | ... | ... | ... | ... | ... |

In place of each composition $x_i * x_j$, there is that element from $M$ which is obtained as a result.

*Example 16.* If the elements of group $S_3$ are designated

$\alpha_0 = e = \begin{pmatrix} 1 & 2 & 3 \\ 1 & 2 & 3 \end{pmatrix}, \alpha_1 = \begin{pmatrix} 1 & 2 & 3 \\ 1 & 3 & 2 \end{pmatrix}, \alpha_2 = \begin{pmatrix} 1 & 2 & 3 \\ 2 & 1 & 3 \end{pmatrix}, \alpha_3 = \begin{pmatrix} 1 & 2 & 3 \\ 2 & 3 & 1 \end{pmatrix}$,

$\alpha_4 = \begin{pmatrix} 1 & 2 & 3 \\ 3 & 1 & 2 \end{pmatrix}, \alpha_5 = \begin{pmatrix} 1 & 2 & 3 \\ 3 & 2 & 1 \end{pmatrix}$,

then table 2 is the Cayley table of $S_3$.

*Remark 8.* A Cayley table describes the structure of a finite group by arranging all the possible products of all its elements in a square table reminiscent of a multiplication table. Many properties of a group - such as whether or not it is abelian, which elements are inverses of which elements, and some other properties - can be discovered from its Cayley table.

Table 2

| · | e | $\alpha_1$ | $\alpha_2$ | $\alpha_3$ | $\alpha_4$ | $\alpha_5$ |
|---|---|---|---|---|---|---|
| e | e | $\alpha_1$ | $\alpha_2$ | $\alpha_3$ | $\alpha_4$ | $\alpha_5$ |
| $\alpha_1$ | $\alpha_1$ | e | $\alpha_3$ | $\alpha_2$ | $\alpha_5$ | $\alpha_4$ |
| $\alpha_2$ | $\alpha_2$ | $\alpha_4$ | e | $\alpha_5$ | $\alpha_1$ | $\alpha_3$ |
| $\alpha_3$ | $\alpha_3$ | $\alpha_5$ | $\alpha_1$ | $\alpha_4$ | e | $\alpha_2$ |
| $\alpha_4$ | $\alpha_4$ | $\alpha_2$ | $\alpha_5$ | e | $\alpha_3$ | $\alpha_1$ |
| $\alpha_5$ | $\alpha_5$ | $\alpha_3$ | $\alpha_4$ | $\alpha_1$ | $\alpha_2$ | e |

A group defined only by the law of composition is called *abstract*. Not every square table, whose cells (except the upper left) are somehow filled with various elements, serves as a Cayley table of some abstract group $(M, *)$: the filling must be such that all four conditions of the definition of a group are fulfilled.

The number of elements $|M|$ of a group $M$ is called its *order*; a group with an infinite set $M$ is considered a *group of infinite order*. For instance, the group $S_n$ of permutations of $n$th degree has order $n!$ and the groups $(\mathbb{Z}, +)$, $(\mathbb{Q}, +)$, $(\mathbb{Q} \setminus \{0\}, \cdot)$ are of infinite order.

SUBSYSTEMS AND ISOMORPHISM OF ALGEBRAIC SYSTEMS. Let a set $M \neq \emptyset$ form an algebraic system with respect to a set of some operations. A non-empty subset $M' \subseteq M$ is called a *subsystem* of this system if, with respect to the same operations, it forms a system of the same type, where each operation on $M'$ is *induced* by the same operation on the whole $M$. The meaning of these two terms (subsystem and induced) will be clarified for one or another algebraic system separately, but the second of them is well illustrated by the following example from mathematical analysis:

*Example 17.* (a) $f_1(x) = \sin x$ on $(-\infty, +\infty)$  (b) $f_2(x) = \sin x$ on $\left[-\frac{\pi}{2}, \frac{\pi}{2}\right]$

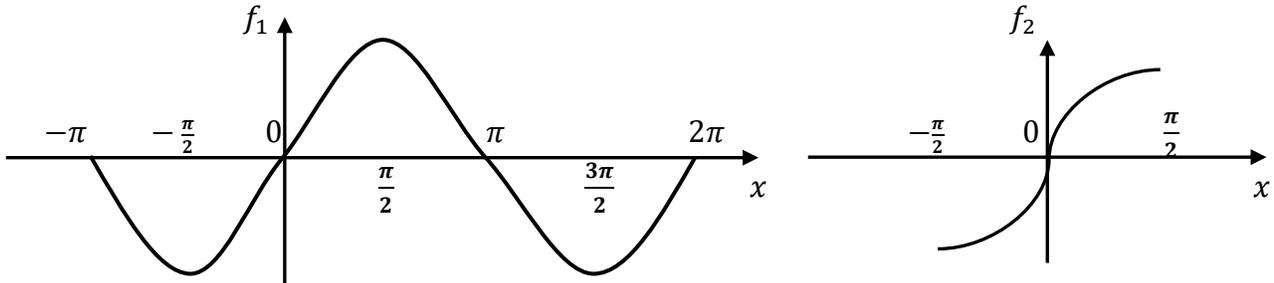

Fig. 7

The functions $f_1$ and $f_2$ (Fig. 7 (a) and (b)) are two different functions since their domains are different. The function $f_2(x)$ is induced by the function $f_1(x)$ on a subset $\left[-\frac{\pi}{2}, \frac{\pi}{2}\right] \subset (-\infty, +\infty)$. This means that for $x \in \left[-\frac{\pi}{2}, \frac{\pi}{2}\right]$, always $f_1(x) = f_2(x)$.

*Remark 9.* Because the sine function is not one-to-one, that is, it is not an injective mapping, we must restrict its domain if, for instance, we want to define its inverse. $f(x) = \sin x$ is one-to-one on $\left[-\frac{\pi}{2}, \frac{\pi}{2}\right]$. The inverse is the arcsine function.

If a binary operation $*$ is defined on a set $M$ and a binary operation $*'$ is defined on its subset $M' \subset M$, then the latter is induced by the operation $*$ in the case when $x *' y = x * y$ for any $x, y \in M'$, in other words, for elements from subset $M'$, the operation $*'$ leads to the same result as $*$ in $M$. In this case, the prime of $*'$ can be omitted, that is, we can denote the operation induced on the subset by the same sign $*$ as in the whole set (similar to the $\sin x$ notation for $f_1(x)$ and $f_2(x)$ described in example 17).

Let $(M, *)$ be a groupoid, $M' \subset M$, $M' \neq \emptyset$. Then $(M', *)$ is a *subgroupoid* of this groupoid if

1. $x, y \in M' \Rightarrow x * y \in M'$ (implying that for any $x, y \in M'$, the quantifier prefix $\forall$ in front is often omitted), that is, $(M', *)$ is also a groupoid;
2. composition $*$ on $M'$ is induced by the composition on $M$.

Subsemigroup and submonoid are defined similarly.

Now, let $(M, *)$ be a group, $M' \subset M$, $M' \neq \emptyset$, and the operation $*$ on $M'$ be induced. To establish that $(M', *)$ is a subgroup of group $(M, *)$, it is sufficient to check in $M'$ not all four conditions of the definition of a group, but only two; 1) and 4). Indeed:
1), that is, the condition $x, y \in M' \Rightarrow x * y \in M'$ is necessary to check for all possible $x, y \in M'$, otherwise $(M', *)$ may not even be a groupoid;
2) (associativity) in $M'$ is performed automatically since it holds for any $x, y, z \in M$, and the operation $*$ on $M'$ is induced;
3) follows from 4): let's take any $x \in M' \neq \emptyset$. According to 4), $x^{-1} \in M'$, and then, due to 1), $e = x * x^{-1} \in M'$.

*Examples 18*. In $S_n$, let's consider the set $S_n(1)$ of all permutations that leave 1 fixed in its place, that is, having the form $\begin{pmatrix} 1 & 2 & \dots & n \\ 1 & i_2 & \dots & i_n \end{pmatrix}$. It is clear that
1. $\alpha, \beta \in S_n(1) \Rightarrow \alpha \beta \in S_n(1)$;
4. if $\alpha = \begin{pmatrix} 1 & 2 & \dots & n \\ 1 & i_2 & \dots & i_n \end{pmatrix} \in S_n(1)$, then $\alpha^{-1} = \begin{pmatrix} 1 & i_2 & \dots & i_n \\ 1 & 2 & \dots & n \end{pmatrix} \in S_n(1)$.

Therefore, $S_n(1)$ is a subgroup of group $S_n$. The set $S_n(1, 2)$ of all permutations of the form $\begin{pmatrix} 1 & 2 & 3 & \dots & n \\ 1 & 2 & i_3 & \dots & i_n \end{pmatrix}$ that leave two elements (1 and 2) fixed in their place is a subgroup of $S_n(1)$ and, explicitly, of group $S_n$. In general, if we consider some set of specific fixed elements $j_1, j_2, \dots, j_k \in \{1, 2, \dots, n\}$ ($0 \leq k \leq n$), then the subset $S_n(j_1, j_2, \dots, j_k) \subseteq S_n$ of all permutations of $n$th degree that leave each of the numbers $j_1, j_2, \dots, j_k$ unchanged will be a subgroup of $S_n$; extreme cases: when $k = 0$, this subgroup coincides with $S_n$ itself and when $k = n$, it consists of only one element $e = \begin{pmatrix} 1 & 2 & \dots & n \\ 1 & 2 & \dots & n \end{pmatrix}$ (*unitary subgroup*).

Let $(M, *)$ be a monoid, $M'$ be a subset of all its invertible elements; $M' \neq \emptyset$ since $e \in M'$. We will show that $(M', *)$ is a submonoid of monoid $(M, *)$.
1. $(M', *)$ is a subgrupoid of $(M, *)$: $x, y \in M' \Rightarrow x * y \in M'$ since the composition of invertible elements is also invertible;
2. associativity in $(M', *)$ automatically follows from associativity in $(M, *)$.

3. $e \in M'$ (since the neutral element is invertible).

So, $(M', *)$ is a submonoid. Moreover, it is a group since, by the definition of $M'$, all its elements are invertible.

Two algebraic systems are called *isomorphic* if there exists a bijection between their sets and a bijection between their operations that preserves their arities such that the corresponding operations on the corresponding elements lead to the corresponding elements. In other words, in isomorphic systems all operations "act in the same way". In general, the exact definition is cumbersome and we will give it (as well as the definition of a subsystem) for each system separately.

Groupoids $(M, *)$ and $(\widetilde{M}, \widetilde{*})$ are *isomorphic* if there exists a bijection $f: M \to \widetilde{M}$ such that

$$\forall x, y \in M: f(x * y) = f(x) \widetilde{*} f(y).$$

*Simplified Notation*. Although the ~ sign cannot be omitted above $M$, it can be omitted above $*$, because from this, it is clear what elements are connected by the composition sign, that is, in which of the two groupoids this operation is considered:

$$\underset{\in M}{f(x)} \underset{\text{in } M}{*} \underset{\in M}{y)} = \underset{\in \widetilde{M}}{f(x)} \underset{\text{in } \widetilde{M}}{*} \underset{\in \widetilde{M}}{f(y)}.$$

Isomorphism is a symmetric (mutual) relation: if the first groupoid is isomorphic to the second, then the second is isomorphic to the first. Indeed, let $f: M \to \widetilde{M}$ be an isomorphism. Then the mapping $f^{-1}: \widetilde{M} \to M$ that relates to each $\tilde{x} \in \widetilde{M}$ the unique $x \in M$, for which $f(x) = \tilde{x}$, is a bijection (why?), and isomorphism: $f^{-1}(\tilde{x} * \tilde{y}) = f^{-1}(f(x) * f(y)) = f^{-1}(f(x * y)) = x * y = f^{-1}(\tilde{x}) * f^{-1}(\tilde{y})$ for any $\tilde{x}, \tilde{y} \in \widetilde{M}$ and also for those $x, y \in M$, for which $f(x) = \tilde{x}$ and $f(y) = \tilde{y}$.

**Theorem 6.** *Let $(M, *)$ and $(\widetilde{M}, \widetilde{*})$ be isomorphic groupoids. If the first of them has any of the properties: commutativity, associativity, existence of a neutral element, or invertibility of all elements, then the second one also has the same property.*

*Proof. Commutativity*. Let $x * y = y * x$ be always in $M$, let $\tilde{x}, \tilde{y} \in \widetilde{M}$ be arbitrary, and let *f* be an isomorphism of groupoids. For $\tilde{x}$ and $\tilde{y}$ there are such preimages $x, y \in M$ (they are also unique, which is not important in this case), for which $f(x) = \tilde{x}$ and $f(y) = \tilde{y}$. We have
$\tilde{x} * \tilde{y} = f(x) * f(y) = f(x * y) = f(y * x) = f(y) * f(x) = \tilde{y} * \tilde{x}$ (Explain each step!).

*Associativity*. Let $x * (y * z) = (x * y) * z$ be always in $M$, let $\tilde{x}, \tilde{y}, \tilde{z} \in \widetilde{M}$ be arbitrary, and let $x, y, z \in M$ be their preimages under the isomorphism $f$ of groupoids. Then $\tilde{x} * (\tilde{y} * \tilde{z}) =$
$= f(x) * (f(y) * f(z)) = f(x) * f(y * z) = f(x * (y * z)) = f((x * y) * z) =$
$= f(x * y) * f(z) = (f(x) * f(y)) * f(z) = (\tilde{x} * \tilde{y}) * \tilde{z}.$

*Neutral element.* Let $e$ be a neutral element of the groupoid $(M, *)$. We will show that $\tilde{e} = f(e)$ is a neutral element in $(\widetilde{M}, *)$. Let's take $\tilde{x} \in \widetilde{M}$ and let $x$ be its preimage under the isomorphism $f$. The equalities

$$\tilde{x} * \tilde{e} = f(x) * f(e) = f(x * e) = f(x) = \tilde{x},$$
$$\tilde{e} * \tilde{x} = f(e) * f(x) = f(e * x) = f(x) = \tilde{x}$$

are true for any $\tilde{x} \in \widetilde{M}$. Therefore $\tilde{e}$ satisfies the definition of a neutral element of a groupoid $(\widetilde{M}, *)$.

*Invertibility.* Let a groupoid $(M, *)$ have a neutral element $e$. Then the isomorphic groupoid $(\widetilde{M}, *)$ also has such an element $\tilde{e} = f(e)$ and the definition of the inverse element makes sense for both groupoids. If $x \in M$ is invertible, then $f(x) \in \widetilde{M}$ is also invertible and the inverse element $[f(x)]^{-1}$ for it is $f(x^{-1})$. In fact, from $x * x^{-1} = x^{-1} * x = e$, it follows that $f(x) * f(x^{-1}) = f(x * x^{-1}) = f(e) = \tilde{e}$ and $f(x^{-1}) * f(x) = f(x^{-1} * x) = f(e) = \tilde{e}$, that is, $f(x^{-1})$ satisfies the definition of the inverse element for $f(x) \in \widetilde{M}$. Therefore, if all elements in $(M, *)$ are invertible, then so are all elements in $(\widetilde{M}, *)$. □

**Corollary 1.** *If one of two isomorphic groupoids is a group* (*an abelian group*), *then the second is also a group* (*respectively, an abelian group*). □

*Example 19*. Let B be a set of six motions of an equilateral triangle $ABC$ (Fig. 8) that transform it into itself: $\beta_0$ is at rest; $\beta_1, \beta_2, \beta_3$ are mirror images (reflections) about the bisectors of angles $A$, $B$, $C$; $\beta_4$ is a rotation around the center $120°$ clockwise, $\beta_5$ is a similar rotation counterclockwise.

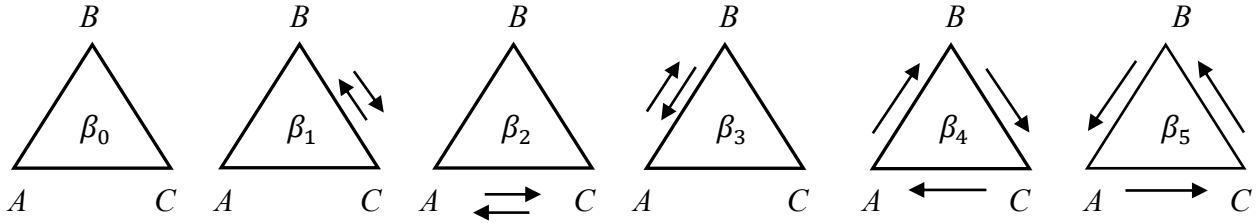

Fig. 8

We call the product $\beta_i \cdot \beta_j$ a movement consisting of the sequential execution of $\beta_i$ and $\beta_j$. For $\beta_0, \beta_1, \beta_2, \beta_3, \beta_4, \beta_5$, we directly compile a multiplication table (table 3).

Table 3

| · | $\beta_0$ | $\beta_1$ | $\beta_2$ | $\beta_3$ | $\beta_4$ | $\beta_5$ |
|---|---|---|---|---|---|---|
| $\beta_0$ | $\beta_0$ | $\beta_1$ | $\beta_2$ | $\beta_3$ | $\beta_4$ | $\beta_5$ |
| $\beta_1$ | $\beta_1$ | $\beta_0$ | $\beta_5$ | $\beta_4$ | $\beta_3$ | $\beta_2$ |
| $\beta_2$ | $\beta_2$ | $\beta_4$ | $\beta_0$ | $\beta_5$ | $\beta_1$ | $\beta_3$ |
| $\beta_3$ | $\beta_3$ | $\beta_5$ | $\beta_4$ | $\beta_0$ | $\beta_2$ | $\beta_1$ |
| $\beta_4$ | $\beta_4$ | $\beta_2$ | $\beta_3$ | $\beta_1$ | $\beta_5$ | $\beta_0$ |
| $\beta_5$ | $\beta_5$ | $\beta_3$ | $\beta_1$ | $\beta_2$ | $\beta_0$ | $\beta_4$ |

Since all its cells are filled with elements of the set B = $\{\beta_0, \beta_1, \beta_2, \beta_3, \beta_4, \beta_5\}$, (B, ·) is a groupoid. Once the isomorphism $f: S_3 \to B$ is established, we will obtain due to Corollary 1 that (B, ·) is a group. To make it easier to find $f$, we will designate the vertices of triangle $ABC$ with the numbers 1, 2, 3 (Fig. 9). Then it is clear that we need to set

$$f(\alpha_0) = \beta_0,\ f(\alpha_1) = \beta_1,\ f(\alpha_2) = \beta_3,\ f(\alpha_3) = \beta_4,\ f(\alpha_4) = \beta_5,\ f(\alpha_5) = \beta_2.$$

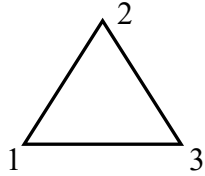

Fig. 9

Consequently (see Fig. 8),

$$f(\alpha_0) = \beta_0 = \begin{pmatrix} A & B & C \\ A & B & C \end{pmatrix},\ f(\alpha_1) = \beta_1 = \begin{pmatrix} A & B & C \\ A & C & B \end{pmatrix},\ f(\alpha_2) = \beta_3 = \begin{pmatrix} A & B & C \\ B & A & C \end{pmatrix},$$

$$f(\alpha_3) = \beta_4 = \begin{pmatrix} A & B & C \\ C & A & B \end{pmatrix},\ f(\alpha_4) = \beta_5 = \begin{pmatrix} A & B & C \\ B & C & A \end{pmatrix},\ f(\alpha_5) = \beta_2 = \begin{pmatrix} A & B & C \\ C & B & A \end{pmatrix}.$$

It is easy to convince yourself that when making $f(\alpha_i \cdot \alpha_j) = f(\alpha_i) \cdot f(\alpha_j)$ for any $\alpha_i, \alpha_j$, $i = 0, \ldots, 5$, the $\beta_i$ values to the right and below the lines of table 3 coincide with those ones of the Cayley table of group $S_3$.

Indeed, let's consider, for example, $f(\alpha_3 \cdot \alpha_3) = f(\alpha_3) \cdot f(\alpha_3) = \beta_4 \cdot \beta_4 =$
$= \begin{pmatrix} A & B & C \\ C & A & B \end{pmatrix} \cdot \begin{pmatrix} A & B & C \\ C & A & B \end{pmatrix} = \begin{pmatrix} A & B & C \\ C & A & B \end{pmatrix} \cdot \begin{pmatrix} C & A & B \\ B & C & A \end{pmatrix} = \begin{pmatrix} A & B & C \\ B & C & A \end{pmatrix} = \beta_5,$
that is, $\beta_4 \cdot \beta_4 = \beta_5$.

On the other hand, $f(\alpha_5 \cdot \alpha_4) = f(\alpha_5) \cdot f(\alpha_4) = \beta_2 \cdot \beta_5 = \begin{pmatrix} A & B & C \\ C & B & A \end{pmatrix} \cdot \begin{pmatrix} A & B & C \\ B & C & A \end{pmatrix} =$
$= \begin{pmatrix} A & B & C \\ C & B & A \end{pmatrix} \cdot \begin{pmatrix} C & B & A \\ A & C & B \end{pmatrix} = \begin{pmatrix} A & B & C \\ A & C & B \end{pmatrix} = \beta_1,$ that is, $\beta_2 \cdot \beta_5 = \beta_1$. Recall that $S_3$ is non-commutative. So, $\beta_2$ is read above the horizontal line of table 3 while $\beta_5$ – to the left of its vertical line. At last, $f(\alpha_4 \cdot \alpha_5) = f(\alpha_4) \cdot f(\alpha_5) = \beta_5 \cdot \beta_2 = \begin{pmatrix} A & B & C \\ B & C & A \end{pmatrix} \cdot \begin{pmatrix} A & B & C \\ C & B & A \end{pmatrix} =$
$= \begin{pmatrix} A & B & C \\ B & C & A \end{pmatrix} \cdot \begin{pmatrix} B & C & A \\ B & A & C \end{pmatrix} = \begin{pmatrix} A & B & C \\ B & A & C \end{pmatrix} = \beta_3,$ that is, $\beta_5 \cdot \beta_2 = \beta_3$.

And so, $(S_3, \cdot)$ and $(B, \cdot)$ are isomorphic. Therefore, (B, ·) is a group.

*Remark 10*. Rearranging the rows and columns in the multiplication table of groupoid (B, ·) so as to give the outer elements (to the left and above the lines of table 3) the order $\beta_0, \beta_1, \beta_3, \beta_4, \beta_5, \beta_2$ and comparing the resulting table with the Cayley table of group $S_3$, we will be convinced that $f(\alpha_i \cdot \alpha_j) = f(\alpha_i) \cdot f(\alpha_j)$ for any $\alpha_i, \alpha_j \in S_3$.

## EXERCISE

Let $F = \{f_1, f_2, f_3, f_4, f_5, f_6\}$ be a set of six functions: $f_1(x) = x$, $f_2(x) = \frac{x-1}{x}$, $f_3(x) = \frac{1}{1-x}$, $f_4(x) = \frac{1}{x}$, $f_5(x) = 1 - x$, $f_6(x) = \frac{x}{x-1}$, and $*$ be a superposition of functions (permutation of functions into a function): $(f_i * f_j)(x) = f_j(f_i(x))$. Show that $(F, *)$ is a groupoid isomorphic to $S_3$, that is, it is a group.

## ALGEBRAIC SYSTEMS WITH TWO BINARY OPERATIONS

Let two binary operations be defined on a set $M \neq \emptyset$, the first of which, we will call *addition* and write *additively*, and the second of which, we will call *multiplication* and write *multiplicatively* despite the fact that they do not necessarily coincide with the ordinary arithmetic operations (and in general $M$ is not necessarily a numerical set). An algebraic system $(M, +, \cdot)$ is called a *ring* if:

I. $(M, +)$ is an abelian group;
II. $(M, \cdot)$ is a groupoid;
III. the operations of addition and multiplication are related by two distributive laws (or distributive properties):

$$\left.\begin{array}{l} x(y + z) = xy + xz, \\ (y + z)x = yx + zx \end{array}\right\} \text{ for any } x, y, z \in M.$$

A ring $(M, +, \cdot)$ is called *associative* if the multiplication operation is associative, that is, $(M, \cdot)$ is a semigroup.

A ring $(M, +, \cdot)$ is called *commutative* if the multiplication operation is commutative. In this case, of the two distributive properties, it is sufficient to check only one since the other will follow from it. Let, for example, it be given that $x(y + z) = xy + xz$ for any $x, y, z \in M$, then $(y + z)x = x(y + z) = xy + xz = yx + zx$, and exactly in the same way the first property is obtained from the second.

In the presence of the multiplication neutral element 1, $(M, +, \cdot)$ itself is called a *ring with unity*.

*Examples 20.* With respect to ordinary arithmetic operations, $(\mathbb{N}, +, \cdot)$ is not a ring since $(\mathbb{N}, +)$ is not a group. Note that conditions II and III are satisfied here ("one wound is fatal, but the other two, fortunately, do not cause concern"). Even the associativity and commutativity of multiplication are evident as well as the existence of a unit.

$(\mathbb{Z}, +, \cdot)$ is an associative-commutative ring with unity. Replacing $\mathbb{Z}$ here with the set of all even numbers only, we get an example of an associative-commutative ring without unity (and what happens in the case of the set of all odd numbers?).
$(\mathbb{Q}, +, \cdot)$ is an associative-commutative ring with unity.
$(M, +, \times)$, where $M$ is the set of all geometric vectors of the three-dimensional Euclidean space, $+$ is the vector addition, $\times$ is the cross product (or vector product), is a ring. As follows from the well-known rules of vector algebra, $(M, +)$ is an abelian group while $(M, \times)$ is neither commutative nor associative, and has no unit either. Instead of commutativity, it has *anticommutativity* $X \times Y = -(Y \times X)$, due to which here also both distributive properties follow from each other.

In non-associative rings "almost nothing can be done" (we are talking about algebraic transformations of some expressions into others) if the lack of associativity is not compensated for by the presence of some other laws. In the example with vectors, such a "compensator" is the Jacobi identity

$$X \times (Y \times Z) + Y \times (Z \times X) + Z \times (X \times Y) = (X \times Y) \times Z + (Y \times Z) \times X + (Z \times X) \times Y = 0.$$

(A ring $(M, +, \cdot)$ is called a *Lie ring* if $xx = 0$ and $(xy)z + (yz)x + (zx)y = 0$ for any $x, y, z \in M$. We suggest that you independently prove that the relation $xx = 0$ implies the anticommutativity).

(F, +, ·), where F is the set of all possible functions $f : \mathbb{R} \to \mathbb{R}$ defined on the real axis and taking real values, and for any $f_1, f_2 \in F$ the function $f_1 + f_2$ assigns to each $x \in \mathbb{R}$ the value $f_1(x) + f_2(x)$, and the function $f_1 f_2$ - the value $f_1(x) \cdot f_2(x)$, represents an associative-commutative ring with unity. The role of the latter is obviously played by the function $f(x) \equiv 1$, which takes the value 1 for all $x \in \mathbb{R}$ (and the role of zero is played by the function $f(x) \equiv 0$). The laws of commutativity, associativity and distributivity of addition and multiplication of functions follow directly from the same laws for real numbers.

Both examples (with vectors and functions) illustrate a possibility that does not occur for number rings (with the usual operations of addition and multiplication). In general, there can be two elements, both nonzero, whose product is zero. So, $X \times X = 0$ even if $X$ is a non-zero vector, and the two functions

$$f_1(x) = \frac{1}{2}(|x| + x) = \begin{cases} x \text{ if } x \geq 0, \\ 0 \text{ if } x \leq 0 \end{cases}$$

and

$$f_2(x) = \frac{1}{2}(|x| - x) = \begin{cases} 0 \text{ if } x \geq 0, \\ -x \text{ if } x \leq 0 \end{cases}$$

of the ring (F, +, ·) (both different from zero) give zero as their product (see Fig. 10)

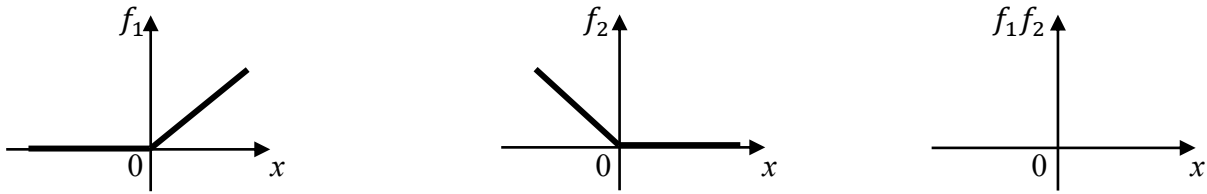

Fig. 10

Elements $x$ and $y$ of a ring $(M, +, \cdot)$ are called *zero divisors* if $x \neq 0, y \neq 0$, but $xy = 0$, where $x$ is a *left* and $y$ is a *right zero divisor* (for commutative rings this distinction obviously loses its meaning). In the previous examples, $(\mathbb{Z}, +, \cdot)$ and $(\mathbb{Q}, +, \cdot)$ are rings without zero divisors, while $(M, +, \times)$ and $(F, +, \cdot)$ are rings with zero divisors.

Let's now proceed to deduce the simplest properties of an arbitrary ring $(M, +, \cdot)$ from its general definition.

(1) $x \cdot 0 = 0 \cdot x = 0$ for any $x \in M$.

Indeed, $x \cdot x + 0 \cdot x = (x + 0) \cdot x = x \cdot x$. Where adding to the first and last parts of the equality the element $-(x \cdot x)$ and using the associativity of addition, we obtain $0 \cdot x = -(x \cdot x) + x \cdot x = = 0$. The case $x \cdot 0 = 0$ can be proved similarly.

(2) $(-x)y = x(-y) = -xy$ for any $x, y \in M$.

The relation $(-x)y = -xy$ means that the element $(-x)y$ is the opposite of $xy$, and this is easy to verify: $xy + (-x)y = [x + (-x)]y = 0 \cdot y = 0$ due to (1). Since $(M, +)$ is an abelian group, $(-x)y + xy = 0$. Equalities $xy + x(-y) = x(-y) + xy = 0$ can be proved similarly.

(3) $(-x)(-y) = xy$ for any $x, y \in M$.

Indeed, applying (2) twice and taking into account that the opposite element of the opposite is the original element, we have $(-x)(-y) = -x(-y) = -(-xy) = xy$.

(4) $\left(\sum_{i=1}^{m} x_i\right)\left(\sum_{j=1}^{n} y_j\right) = \sum_{i=1}^{m} \sum_{j=1}^{n} x_i y_j$ for any $x_1, x_2, \ldots, x_m, y_1, y_2, \ldots, y_n \in M; m, n \in \mathbb{N}$.

This *general law of distributivity* is proven by $(mn - 1)$-multiple applications of the laws III (using the associativity of addition):

$$(x_1 + x_2 + \cdots + x_m)(y_1 + y_2 + \cdots + y_n) = [x_1 + (x_2 + \cdots + x_m)](y_1 + y_2 + \cdots + y_n) =$$
$$= x_1(y_1 + y_2 + \cdots + y_n) + (x_2 + \cdots + x_m)(y_1 + y_2 + \cdots + y_n) = \cdots =$$
$$= x_1(y_1 + y_2 + \cdots + y_n) + x_2(y_1 + y_2 + \cdots + y_n) + \cdots + x_m(y_1 + y_2 + \cdots + y_n) = \cdots =$$
$$= x_1 y_1 + x_1(y_2 + \cdots + y_n) + \cdots = \cdots =$$
$$= x_1 y_1 + x_1 y_2 + \cdots + x_1 y_n + x_2 y_1 + x_2 y_2 + \cdots + x_2 y_n + \cdots + x_m y_1 + x_m y_2 + \cdots + x_m y_n.$$

*Remark 11.* Properties (2) and (3) express the "rule of signs in multiplication", but it is necessary to remember that $-x$ means the element opposite to $x$ (that is, the inverse of addition), and not at all "negative": even for ordinary numbers the number $-x$ is not necessarily negative (after all, it can be that $x < 0$). In a ring of general type, the concept of "positivity" and "negativity" of an element has no meaning. The sum $x + (-y)$ can be written simply as $x - y$ and called the *difference* of elements $x$ and $y$.

**Theorem 7.** *In a ring without zero divisors, each of the equalities $cx = cy$ and $xc = yc$ for $c \neq 0$ implies $x = y$. In other words, the equality can be reduced (from the same side) by a common factor different from zero.* □

*Remark 12.* Indeed, if for example $cx = cy$, then $cx - cy = 0 \Rightarrow c(x - y) = 0$ and since $c \neq 0$ and there are no zero divisors in the ring, it is a must that $x - y = 0$, that is, $x = y$. The case $xc = yc$ is similar.

In a ring with unity, an invertible element (that is, one having an inverse under multiplication) is also called a *unit divisor*.

**Theorem 8.** *In an associative ring $(M, +, \cdot)$ with unity, no zero divisor can be a unit divisor.*

*Proof.* Let's assume, for example, that $xy = 0$, where $y \neq 0$, and $x$ has an inverse element $x^{-1}$ (the case where $x \neq 0$, and $y$ has $y^{-1}$, is considered analogously). Then $x^{-1}(xy) = x^{-1} \cdot 0 = 0$, and because of associativity, $(x^{-1}x)y = 0 \Rightarrow 1 \cdot y = 0 \Rightarrow y = 0$, contrary to our assumption. □

A ring consisting of a single element (which obviously serves as both zero and one) is called *trivial*. If in a ring $(M, +, \cdot)$ there is $1 = 0$, then it is trivial because $x = x \cdot 1 = x \cdot 0 = 0$ for any $x \in M$. A non-trivial associative-commutative ring with unity ($1 \neq 0$), having no zero divisors, is called an *integral ring*. Such, for example, is the ring of integers $(\mathbb{Z}, +, \cdot)$.

An associative ring with unity, in which every nonzero element is invertible, is called a *division ring* (or *skew field*), and if the multiplication operation is commutative, then a *field*. In view of the special importance of the latter, we will give its definition, which does not explicitly rely on the concept of a ring.

An algebraic system $(M, +, \cdot)$ is called a field if:
I. $(M, +)$ is an abelian group;
II. $(M, \cdot)$ is a groupoid and $(M \setminus \{0\}, \cdot)$ is an abelian group;
III. addition and multiplication are related by the distributive property (it can be written in either of the two forms).

The equivalence of both definitions of a field is proven without difficulty (we suggest you do it yourself).

*Examples 21.* $(\mathbb{Z}, +, \cdot)$ is not a field (not even a skew field) since integers, including 0, do not have inverses under multiplication. However, the only unit divisors here are 1 and -1 (prove it!).
$(\mathbb{Q}, +, \cdot)$ is a field.

Due to Theorem 8 on zero and unit divisors, *no skew field*, in particular, *no field has zero divisors*. Therefore, the rings of $(M, +, \times)$ vectors and $(F, +, \cdot)$ functions are not skew fields, much less fields. A field cannot be a trivial ring since $1 \in M \setminus \{0\}$ and, therefore, $1 \neq 0$.

INTERPRETATION OF THE SIMPLEST PROPERTIES OF A FIELD FROM ITS DEFINITION. Since the multiplication operation in a field is commutative, $xy^{-1} = y^{-1}x$ for any $y \neq 0$ and writing each of these products as a fraction $\frac{x}{y}$ will not lead to confusion. Below, it is assumed that $x, y, z, t$ are arbitrary elements of a field $(M, +, \cdot)$, and those of them that are placed in the denominator are different from zero.

(1) $\frac{x}{y} = \frac{z}{t}$ if and only if $xt = yz$.

Indeed, from $xy^{-1} = zt^{-1}$ multiplying both parts by $yt$ and using associativity and commutativity, we get $xt = yz$. Conversely, from the last equality, assuming $y \neq 0$ and $t \neq 0$, multiplying by $y^{-1}t^{-1}$, we obtain $xy^{-1} = zt^{-1}$.

(2) $\frac{x}{y} \pm \frac{z}{t} = \frac{xt \pm yz}{yt}$.

In fact, the left-hand side of $xy^{-1} \pm zt^{-1} = xy^{-1}tt^{-1} \pm zt^{-1}yy^{-1} = xt(yt)^{-1} \pm yz(yt)^{-1} = (xt \pm yz)(yt)^{-1}$, that is, is equal to the right-hand side, and from $y \neq 0$ and $t \neq 0$, it follows $yt \neq 0$ due to the absence of a zero divisor, and $(yt)^{-1} = t^{-1}y^{-1}$ due to commutativity.

(3) $\frac{x}{y} \cdot \frac{z}{t} = \frac{xz}{yt}$.

(4) $\left(\frac{x}{y}\right)^{-1} = \frac{x^{-1}}{y^{-1}} = \frac{y}{x}$.

(5) $\frac{-x}{y} = \frac{x}{-y} = -\frac{x}{y}$.

*Remark 13.* We suggest you prove (3) - (5) yourself. Identities (2) – (5) reveal the "rules for operations with fractions".

THE CHARACTERISTIC OF A FIELD. In numerical fields with the usual arithmetic operations, the natural numbers $1$, $1 + 1 = 2$, $1 + 1 + 1 = 3$, and so on, all are different. But in an arbitrary field $(M, +, \cdot)$, the sums of $n$ and $n + p$ units ($p \geq 1$) may be the same, or equivalently, the sum of $p$ units may be equal to zero (for $p = 1$ this is the "sum of one term"). The smallest number of $p$ units that add up to $0$ is called the *characteristic* of the field. If $1 \neq 0$, $1 + 1 \neq 0$, $1 + 1 + 1 \neq 0$, ... for any number of terms, then the field characteristic is by definition considered equal to zero.

*Examples 22.* Field $(\mathbb{Q}, +, \cdot)$ with the usual addition and multiplication has characteristic $p = 0$.

For $p = 1$, we would have a trivial ring, which, as is known, is not a field. Therefore, *fields of characteristic* 1 *do not exist*.

Algebraic system $(\{0, 1\}, +, \cdot)$ with the operations

$0 + 0 = 0, \qquad 0 + 1 = 1 + 0 = 1, \qquad 1 + 1 = 0, \qquad 0 \cdot 0 = 0 \cdot 1 = 1 \cdot 0 = 0, \qquad 1 \cdot 1 = 1$

(The old Russian arithmetic "even – odd") is a field of characteristic 2. We suggest checking yourself the fulfillment of all points of the field definition. Examples of fields with other characteristic values will be considered later.

**Theorem 9.** *If the characteristic $p$ of a field is nonzero, then $p$ is a prime number.*

*Proof.* Let's assume the opposite: the characteristic of some field has the form $p = mn$, where $m$ and $n$ are natural numbers such that $1 < m < p$ and $1 < n < p$. The sum of $p$ units is $1 + 1 + \cdots + 1 = 0$. The terms of this sum can be partitioned into $m$ groups of $n$ units each. But then

$$\underbrace{1 + 1 + \cdots + 1}_{p} = \underbrace{(1 + 1 + \cdots + 1)}_{n} + \underbrace{(1 + 1 + \cdots + 1)}_{n} + \cdots + \underbrace{(1 + 1 + \cdots + 1)}_{n} =$$
$$= 1 \cdot (1 + 1 + \cdots + 1) + 1 \cdot (1 + 1 + \cdots + 1) + \cdots + 1 \cdot (1 + 1 + \cdots + 1)$$
$$= \underbrace{(1 + 1 + \cdots + 1)}_{m} \cdot \underbrace{(1 + 1 + \cdots + 1)}_{n}$$

based on the distributive property. Since $p$ is the smallest number of units that add up to zero and $m < p, n < p$, $\underbrace{1 + 1 + \cdots + 1}_{m} \neq 0$ and $\underbrace{1 + 1 + \cdots + 1}_{n} \neq 0$, at the same time as $\underbrace{1 + 1 + \cdots + 1}_{p} = 0$. Therefore, substituting both $m \neq 0$ and $n \neq 0$ into $mn = p$ gives as a result $p = 0$, but this is impossible since a field has no zero divisors. □

**Theorem 10.** *If $(M, +, \cdot)$ is a field of characteristic $p$, then $\underbrace{x + x + \cdots + x}_{p} = 0$ for any $x \in M$.* □

*Remark 14.* Indeed, $x + x + \cdots + x = 1 \cdot x + 1 \cdot x + \cdots + 1 \cdot x = \left( \underbrace{1 + 1 + \cdots + 1}_{p} \right) x = 0 \cdot x = 0$ when $p \neq 0$, and "the sum of a zero number of terms" is considered to be equal to zero by definition.

ABOUT RESIDUE RINGS. Let's consider the algebraic system $\mathbb{Z}_m = (\{0, 1, 2, \ldots, m - 1\}, +, \cdot)$, where $m \in \mathbb{N} \setminus \{1\}$, and addition and multiplication are defined as follows: first, the usual arithmetic operation is performed on a pair of numbers, then the result is replaced by its remainder from the division by $m$. For example, in $\mathbb{Z}_7 = (\{0, 1, 2, 3, 4, 5, 6\}, +, \cdot)$, we have $1 + 3 = 4, 3 + 6 = 2, 3 + 4 = 0, 2 \cdot 3 = 6, 2 \cdot 4 = 1$ and so on. $\mathbb{Z}_m$ is called a *ring of residue classes modulo $m$*.

As will be proved in the course on number theory (2nd semester), for $m$ prime, $\mathbb{Z}_m$ is a field of characteristic $p = m$. For $m$ composite, it is only a ring (associative-commutative with unity), but not a field due to the presence of zero divisors (thus, in $\mathbb{Z}_6$, we have $2 \cdot 3 = 0$ although $2 \neq 0$ and $3 \neq 0$). The special case of $\mathbb{Z}_2$ was considered above as one of the examples 22.

SUBRING AND SUBFIELD. The general definition of a subsystem of an algebraic system as applied to rings is deciphered as follows: a non-empty subset $M' \subseteq M$ of elements of a ring $(M, +, \cdot)$ forms its *subring* if the system $(M', +, \cdot)$ is a ring with respect to the operations induced on $M'$ by the previous $+$ and $\cdot$ operations (and denoted in the same way. Allowing for freedom, one can say: "with respect to the same operations"). Replacing the word "ring" everywhere with the word "field", we obtain the definition of a *subfield*.

*Examples 23.* The ring of all even numbers is a subring of the ring of integers $(\mathbb{Z}, +, \cdot)$, but the set of odd numbers is not. The field, and in this connection, the ring $(\mathbb{Q}, +, \cdot)$ contains the subring $(\mathbb{Z}, +, \cdot)$ that is not a subfield. None of the residue rings $\mathbb{Z}_m$ for $m > 1$ is a subring of $(\mathbb{Z}, +, \cdot)$ since although $M' = \{0, 1, 2, \ldots, m - 1\}$ is a subset of $M = \mathbb{Z}$, the addition operation (and in the case of $m > 2$, the multiplication operation) in $(M', +, \cdot)$ is not induced by the operation of the same name in $(M, +, \cdot)$, that is, it is a "different operation".

If $M'$ is a non-empty subset of elements of the ring $(M, +, \cdot)$, then to check that $(M', +, \cdot)$ is a subring, it is sufficient to establish the following:

(a) $(M', +)$ and $(M', \cdot)$ are groupoids, that is, $x + y \in M'$ and $x \cdot y \in M'$ always follow from $x, y \in M'$;

(b) If $x \in M'$, then $-x \in M'$.

Actually, since $M' \subseteq M$ and $(M, +)$ is a group, from the first part of (a) and from (b) it already follows that $(M', +)$ is also a group; the second part of (a) ensures condition II: $(M', \cdot)$ is a groupoid; condition III (distributivity) and commutativity of addition are automatically satisfied on subset $M'$ as long as they are satisfied on all $M$.

In the case where $(M, +, \cdot)$ is a field, to check that $(M', +, \cdot)$ is its subfield, it is sufficient to check, in addition to (a) and (b), that the following condition is met:

(c) If $x \in M' \setminus \{0\}$, then $x^{-1} \in M' \setminus \{0\}$.

Indeed, the invertibility of any element of set $M' \setminus \{0\}$ implies the absence of zero divisors in $M'$ which, together with the second part of (a), leads to the conclusion that $(M', +, \cdot)$ is a groupoid. From here and from (c), it follows that this is a group since $M' \setminus \{0\} \subseteq M \setminus \{0\}$ and $(M \setminus \{0\}, \cdot)$ is a group. Finally, the commutativity of multiplication in $M'$ automatically follows from the commutativity in all of $M$.

ISOMORPHISM OF RINGS AND FIELDS. The general concept of isomorphism of algebraic systems as applied to rings is expressed by the following precise definition: a ring $(M, +, \cdot)$ is *isomorphic* to a ring $(\widetilde{M}, +, \cdot)$ if there exists a bijection $f: M \to \widetilde{M}$ such that for any $x, y \in M$

$$f(x + y) = f(x) + f(y) \text{ and } f(x \cdot y) = f(x) \cdot f(y).$$

Note that the operations $+$ and $\cdot$ even with $M = \widetilde{M}$ are not necessarily the same in both rings although they are designated identically (and with $M \cap \widetilde{M} = \emptyset$ the very concept of "the same operation" generally loses its meaning). We also note that in the definition of isomorphism it is actually unnecessary to consider, beforehand, the second algebraic system as a ring - it is enough to assume that $(\widetilde{M}, +)$ and $(\widetilde{M}, \cdot)$ are groupoids since from Theorem 6 on the transfer of the properties of a groupoid to an isomorphic groupoid, it follows that the system $(\widetilde{M}, +, \cdot)$ together with the ring $(M, +, \cdot)$ satisfy conditions I and II. Distributivity is obtained in the following way: if given any $\tilde{x}, \tilde{y}, \tilde{z} \in \widetilde{M}$ and $x, y, z \in M$ such that $f(x) = \tilde{x}$, $f(y) = \tilde{y}$ and $f(z) = \tilde{z}$, then

$$\tilde{x}(\tilde{y} + \tilde{z}) = f(x) \cdot [f(y) + f(z)] = f(x) \cdot f(y + z) = f(x(y + z)) = f(xy + xz) = f(xy) + f(xz) = f(x) \cdot f(y) + f(x) \cdot f(z) = \tilde{x}\tilde{y} + \tilde{x}\tilde{z},$$

and in a similar way $(\tilde{y} + \tilde{z})\tilde{x} = \tilde{y}\tilde{x} + \tilde{z}\tilde{x}$.

From the same Theorem 6 on the transfer of properties, it follows that any of the properties of associativity, commutativity, existence of a unit and invertibility of nonzero elements, being fulfilled in the original ring, is also fulfilled in the isomorphic one. Therefore, *if $(M, +, \cdot)$ is a field, then its isomorphic image $(\widetilde{M}, +, \cdot)$ is also a field*.

## FIELD EXTENSIONS

A field $\hat{P} = (\hat{M}, +, \cdot)$ is called an *extension* of a field $P = (M, +, \cdot)$ if the latter is a subfield of $\hat{P}$, where $M \subset \hat{M}$ (that is, $M \subseteq \hat{M}$, but $M \neq \hat{M}$), in other words, if $P$ is a *proper subfield of $\hat{P}$*.

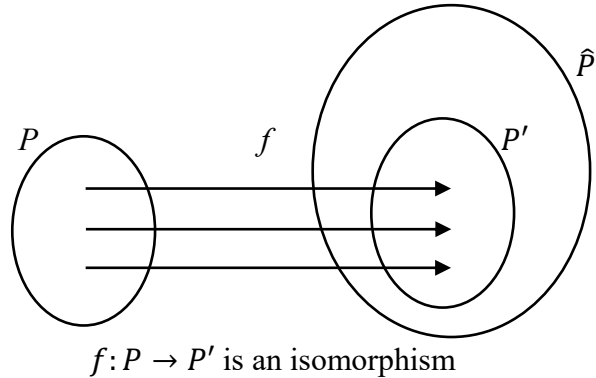

$f: P \to P'$ is an isomorphism

Fig. 11

$P$ is also called $\hat{P}$ if $\hat{P}$ contains not only $P$ itself, but a proper subfield $P'$ isomorphic to $P$ (Fig. 11). Using one particular example, we will demonstrate a typical method of extension in algebra (i.e., the construction of some extension) of a given field.

The historical process of successive extensions of the field $\mathbb{Q}$ of rational numbers consisted in the fact that whenever there arose a need to operate with one or another quantity that could not be expressed by any rational number exactly, but had a clear geometric meaning (the diagonal of a unit square, the length of the circumference of a circle of unit diameter, and others) and allowing this value a rational approximation with an arbitrarily small error; such a value was also called a number although "deaf", "inexpressible", "irrational" (this Latin term can be translated in two ways: "not being a ratio" and "not amenable to understanding"), somehow designated ($\sqrt{2}$, $\pi$ and the like) and acted with it according to the same rules of arithmetic as with previously known numbers. The validity of such extrapolation of laws was confirmed by geometric constructions, as well as approximate calculations "with any predetermined accuracy". So, adding to rational numbers the "dumb" $\sqrt{2}$, for which $\left(\sqrt{2}\right)^2 = 2$, we must include all numbers of the form $x + y\sqrt{2}$, where $x, y \in \mathbb{Q}$, and add and multiply such numbers according to the rules

$$\left.\begin{aligned}(x_1 + y_1\sqrt{2}) + (x_2 + y_2\sqrt{2}) &= (x_1 + x_2) + (y_1 + y_2)\sqrt{2} \\ (x_1 + y_1\sqrt{2}) \cdot (x_2 + y_2\sqrt{2}) &= (x_1 x_2 + 2 y_1 y_2) + (x_1 y_2 + y_1 x_2)\sqrt{2}\end{aligned}\right\} \quad (1)$$

As an exercise we propose to check that the algebraic system $(\{x + y\sqrt{2}/x, y \in \mathbb{Q}\}, +, \cdot)$ is a field (the existence of $(x + y\sqrt{2})^{-1}$ in all cases except $x = y = 0$ can be proved by using the well-known school method of "getting rid of irrationality in the denominator", i.e., "rationalizing the denominator"). In this way, as we will show below, we obtain the "most economical" of the extensions of field $\mathbb{Q}$, leading to the solvability of equation $x^2 = 2$, whose coefficients belong to the field $\mathbb{Q}$, but which has no roots in this field.

The historical path of joining $\sqrt{2}$ to the field $\mathbb{Q}$ cannot be considered erroneous (after all, the results are correct!), but it is not justified: why, for example, the previously unknown number, the square of which is equal to 2, can be added to rational numbers, but a number that is simultaneously equal to 2 and 3 cannot? In the second case, the contradiction is immediately visible, but where is the guarantee that in the first case it will not be discovered someday? All these doubts will disappear if we construct from some clearly defined elements such a field $\widehat{\mathbb{Q}}$ that contains a subfield $\mathbb{Q}'$ isomorphic to $\mathbb{Q}$, but in which the equation $x^2 = 2'$ is already solvable, where $2' \in \mathbb{Q}'$ is an element corresponding to the number $2 \in \mathbb{Q}$ under the isomorphism $f: \mathbb{Q} \to \mathbb{Q}'$.

The fact that a number of the form $x + y\sqrt{2}$ (which arose "spontaneously" in the historical process) is completely determined by assigning an ordered pair $(x, y)$ of rational numbers to it suggests the idea of taking all possible such pairs as elements of $\mathbb{Q}$, and rules (1) suggest (but do not prove!) the following definition of addition and multiplication of pairs:

$$(x_1, y_1) + (x_2, y_2) = (x_1 + x_2, y_1 + y_2),$$
$$(x_1, y_1) \cdot (x_2, y_2) = (x_1 x_2 + 2 y_1 y_2, x_1 y_2 + y_1 x_2).$$

*Algebraic system* $\widehat{\mathbb{Q}} = (\{(x,y)/x, y \in \mathbb{Q}\}, +, \cdot)$ *is a field*. To prove it, we need to check all points (I, II, III) of the definition of a field entirely since pairs of numbers are new objects for us, not belonging to any previously known algebraic system. But the verification of some sub-items of I and II can be combined for shortness. And so, for briefness, let's agree to denote by $\widehat{\mathbb{Q}}$ not only the entire system, but also the set $\{(x, y)/x, y \in \mathbb{Q}\}$ itself.

$(\widehat{\mathbb{Q}}, +)$ and $(\widehat{\mathbb{Q}}, \cdot)$ are *groupoids*. This follows directly from the definition of operations on pairs and from the fact that $(\mathbb{Q}, +)$ and $(\mathbb{Q}, \cdot)$ are groupoids.

*Addition and multiplication in* $\widehat{\mathbb{Q}}$ *are commutative*. Indeed, due to the commutativity of these operations in $\mathbb{Q}$ and their definition in $\widehat{\mathbb{Q}}$,

$$(x_1, y_1) + (x_2, y_2) = (x_1 + x_2, y_1 + y_2) = (x_2 + x_1, y_2 + y_1) = (x_2, y_2) + (x_1, y_1),$$
$$(x_1, y_1) \cdot (x_2, y_2) = (x_1 x_2 + 2 y_1 y_2, x_1 y_2 + y_1 x_2) = (x_2 x_1 + 2 y_2 y_1, x_2 y_1 + y_2 x_1) =$$
$$= (x_2, y_2) \cdot (x_1, y_1).$$

*Addition and multiplication in* $\widehat{\mathbb{Q}}$ *are associative*. Indeed,

$$(x_1, y_1) + [(x_2, y_2) + (x_3, y_3)] = (x_1, y_1) + (x_2 + x_3, y_2 + y_3) =$$
$$= ((x_1 + (x_2 + x_3), y_1 + (y_2 + y_3)) = ((x_1 + x_2) + x_3, (y_1 + y_2) + y_3) =$$
$$= (x_1 + x_2, y_1 + y_2) + (x_3, y_3) =$$
$$= [(x_1, y_1) + (x_2, y_2)] + (x_3, y_3),$$

since addition is associative in $\mathbb{Q}$. Further,

$$(x_1, y_1) \cdot (x_2, y_2)(x_3, y_3) = (x_1, y_1) \cdot (x_2 x_3 + 2 y_2 y_3, x_2 y_3 + y_2 x_3) =$$
$$= \underline{(x_1(x_2 x_3 + 2 y_2 y_3) + 2 y_1 (x_2 y_3 + y_2 x_3), x_1(x_2 y_3 + y_2 y_3) + y_1(x_2 x_3 + 2 y_2 y_3))},$$

$$(x_1, y_1)(x_2, y_2) \cdot (x_3, y_3) = (x_1 x_2 + 2 y_1 y_2, x_1 y_2 + y_1 x_2) \cdot (x_3, y_3) =$$
$$= \underline{((x_1 x_2 + 2 y_1 y_2) x_3 + 2(x_1 y_2 + y_1 x_2) y_3, (x_1 x_2 + 2 y_1 y_2) y_3 + (x_1 y_2 + y_1 x_2) x_3)}$$

(instead of using square brackets, we used the option of omitting or not omitting the dot when denoting the multiplication of pairs). The equality of both resulting pairs underlined by a dashed line, follows from the laws of commutativity, associativity and distributivity in the field $\mathbb{Q}$: we suggest that you do the "school" calculations yourself to prove the equality of both the first and second components of these pairs.

*In $\widehat{\mathbb{Q}}$ the distributive property* III *is fair*, namely,

$$(x_1, y_1)[(x_2, y_2) + (x_3, y_3)] = (x_1, y_1)(x_2 + x_3, y_2 + y_3) =$$
$$= \underline{(x_1(x_2 + x_3) + 2 y_1 (y_2 + y_3), x_1(y_2 + y_3) + y_1(x_2 + x_3))},$$

$$(x_1, y_1)(x_2, y_2) + (x_1, y_1)(x_3, y_3) =$$
$$= (x_1 x_2 + 2 y_1 y_2, x_1 y_2 + y_1 x_2) + (x_1 x_3 + 2 y_1 y_3, x_1 y_3 + y_1 x_3) =$$
$$= \underline{(x_1 x_2 + 2 y_1 y_2 + x_1 x_3 + 2 y_1 y_3, x_1 y_2 + y_1 x_2 + x_1 y_3 + y_1 x_3)},$$

and again, using the calculations in field $\mathbb{Q}$, we are convinced that both obtained pairs (underlined by a dashed line) coincide.

$(\widehat{\mathbb{Q}}, +)$ *is an abelian group*. Points 1 and 2 of the definition of a group and its commutativity have already been proved. But, obviously, the pair $(0,0)$ is the neutral element with respect to addition and the opposite for the pair $(x, y)$ is $(-x, -y)$.

$(\widehat{\mathbb{Q}} \setminus \{(0,0)\}, \cdot)$ *is an abelian group*. From all that has been said above, it already follows that $(\widehat{\mathbb{Q}}, +, \cdot)$ is an associative-commutative ring (check yourself all the points of its definition!). The pair $(1, 0)$ serves as a neutral element under multiplication:

$$(x, y)(1, 0) = (1, 0)(x, y) = (1 \cdot x + 2 \cdot 0 \cdot y, 1 \cdot y + 0 \cdot x) = (x, y).$$

Let's prove the invertibility of each non-zero element. Let $(a, b) \neq (0, 0)$, that is, at least one of the numbers $a, b \in \mathbb{Q}$ is different from zero; we find a pair $(x, y)$, for which $(x, y) \cdot (a, b) = (a, b) \cdot (x, y) = (1, 0)$.

In view of the commutativity of multiplication in $\widehat{\mathbb{Q}}$, it is sufficient to find a pair $(x, y)$ satisfying the equation $(a, b)(x, y) = (1, 0)$; this will be the desired element $(a, b)^{-1}$ since in the monoid $(\widehat{\mathbb{Q}}, \cdot)$ there can only be one of it. To find the pair $(x, y)$, we write the equation that it satisfies in the form
$$(ax + 2by, ay + bx) = (1, 0),$$
which is equivalent to a system of two equations with two unknowns $x$ and $y$ in $\mathbb{Q}$:
$$ax + 2by = 1,$$
$$bx + ay = 0.$$
eliminating $y$ from these equations (by multiplying the first equation by $a$, the second by $-2b$ and adding), we get
$$(a^2 - 2b^2)x = a,$$
and eliminating $x$ from the same equations gives
$$(a^2 - 2b^2)y = -b.$$

The factor $a^2 - 2b^2 \neq 0$, indeed, from $a^2 = 2b^2$ due to the assumption $(a, b) \neq (0, 0)$ it would follow that $b \neq 0$ and, further, $(a/b)^2 = 2$, but this is impossible: there is no rational number $a/b \in \mathbb{Q}$ whose square would be equal to two. Therefore, the desired $x$ and $y$ can only be as follows:
$$x = \frac{a}{a^2 - 2b^2}, \quad y = \frac{-b}{a^2 - 2b^2},$$
and the fact that the resulting pair $(x, y)$ is actually $(a, b)^{-1}$ can be verified directly:
$$(a, b)\left(\frac{a}{a^2 - 2b^2}, \frac{-b}{a^2 - 2b^2}\right) = \left(a \cdot \frac{a}{a^2 - 2b^2} + 2b \cdot \frac{-b}{a^2 - 2b^2}, a \cdot \frac{-b}{a^2 - 2b^2} + b \cdot \frac{a}{a^2 - 2b^2}\right)$$
$$= \left(\frac{a^2 - 2b^2}{a^2 - 2b^2}, \frac{-ab + ab}{a^2 - 2b^2}\right) = (1, 0).$$

Note that for a purely formal proof of the invertibility of the pair $(a, b) \neq (0, 0)$, it would be enough just to carry out this check, taking a pair $\left(\frac{a}{a^2 - 2b^2}, \frac{-b}{a^2 - 2b^2}\right)$ "out of thin air" (and showing that the denominator is different from zero), but we preferred a targeted search to guesswork.

From the invertibility of non-zero pairs, it now follows that $(\mathbb{Q} \setminus \{(0, 0)\}, \cdot)$ is a groupoid: if $(x_1, y_1) \neq (0, 0)$ and $(x_2, y_2) \neq (0, 0)$, then both pairs are invertible and therefore cannot be zero divisors, from where $(x_1, y_1)(x_2, y_2) \neq (0, 0)$.

Thus, all conditions included in points I, II, III of the field definition are verified for the $(\widehat{\mathbb{Q}}, +, \cdot)$ system (only in a different sequence).

*The field of pairs $\widehat{\mathbb{Q}}$ contains a subfield isomorphic to $\mathbb{Q}$.* Let's consider a proper subset $\mathbb{Q}' \subset \widehat{\mathbb{Q}}$ of pairs of the form $(x, 0)$, where $x \in \mathbb{Q}$, and show that it forms a subfield in $\widehat{\mathbb{Q}}$. This time it is sufficient to check three points.

(a) $(\mathbb{Q}', +)$ and $(\mathbb{Q}', \cdot)$ are *groupoids*. Indeed, $(x, 0) + (y, 0) = (x + y, 0) \in \mathbb{Q}'$ and $(x, 0) \cdot (y, 0) = (xy + 2 \cdot 0 \cdot 0, x \cdot 0 + 0 \cdot y) = (xy, 0) \in \mathbb{Q}'$;
(b) $-(x, 0) = (-x, 0) \in \mathbb{Q}'$;
(c) if $(x, 0) \in \mathbb{Q}' \setminus \{(0, 0)\}$, then $x \neq 0$, but then $(x, 0)^{-1} = (x^{-1}, 0) \in \mathbb{Q}' \setminus \{(0, 0)\}$.

To prove that the field $\mathbb{Q}$ of rational numbers is isomorphic to the subfield $\mathbb{Q}'$ of the field of pairs $\widehat{\mathbb{Q}}$, we define a mapping $f \colon \mathbb{Q} \to \mathbb{Q}'$, setting $f(x) = (x, 0)$ for any $x \in \mathbb{Q}$. It is clear that this is a bijection (do the rigorous reasoning yourself!), and the fact that $f$ is an isomorphism is proven without difficulty:

$$f(x + y) = (x + y, 0) = (x, 0) + (y, 0) = f(x) + f(y),$$
$$f(x \cdot y) = (xy, 0) = (x, 0) \cdot (y, 0) = f(x) \cdot f(y).$$

*In the field $\widehat{\mathbb{Q}}$, the equation $(x, y)^2 = (2, 0)$ is solvable.* Let's first note that it is precisely for such an equation that we sought solvability in the extension $\widehat{\mathbb{Q}}$ of the field $\mathbb{Q}$: $(x, y)$ is an unknown element of $\widehat{\mathbb{Q}}$, $(x, y)^2$ means $(x, y) \cdot (x, y)$, and the element $(2, 0) \in \widehat{\mathbb{Q}}$ (belonging to $\mathbb{Q}'$) corresponds to the number $2 \in \mathbb{Q}$ under the natural isomorphism $f \colon \mathbb{Q} \to \mathbb{Q}'$.

To prove the existence of the required pair $(x, y)$, a preliminary search such as for $(a, b)^{-1}$ turns out to be superfluous since "historical considerations" immediately suggest one of the solutions $(x, y) = (0, 1)$: after all, $\sqrt{2}$ is a number of the form $x + y\sqrt{2}$, where $x = 0, y = 1$. For complete rigor only a check is needed:

$$(0, 1)^2 = (0, 1)(0, 1) = (0 \cdot 0 + 2 \cdot 1 \cdot 1, 0 \cdot 1 + 1 \cdot 0) = (2, 0).$$

The desired extension of the field of rational numbers has been constructed. We propose as an exercise to check that the equation $(x, y)^2 = (2, 0)$ is also satisfied by the pair $(0, -1)$, and to prove that there are no other solutions besides these two (here we will need a search that reduces to solving a system of two equations in rational numbers).

Since in an algebraic system it is not the nature of elements that is important (and especially not the way they are written), but the rules for operating on them, $\mathbb{Q}'$ and $\mathbb{Q}$ can be considered the same field, and the pair $(x, 0)$ is just another designation for the rational number $x$; with the same success, the number $x$ itself can be considered as another designation for the pair $(x, 0)$. It is natural to consider symbol $\sqrt{2}$ as the designation of the pair $(0, 1)$, and then it is equally natural to write an arbitrary pair $(x, y)$ in the form $x + y\sqrt{2}$ since $(x, y) = (x, 0) + (0, y) = (x, 0) + (y, 0)(0, 1)$ and such a representation is unique: if $(x, y) = (x', 0) + (y', 0) \cdot (0, 1)$, then the sum is equal to $(x', 0) + (y' \cdot 0 + 2 \cdot 0 \cdot 1, y' \cdot 1 + 0 \cdot 0) = (x', 0) + (0, y') = (x', y')$, from where
$$(x, y) = (x', y'),$$

that is, $x = x'$ and $y = y'$ (equality of ordered pairs!).

Thus, the extension of the field $\mathbb{Q}$, leading to the solvability of the equation $x^2 = 2$, can be considered, instead of the field of pairs $\widehat{\mathbb{Q}}$, the field

$$\mathbb{Q}(\sqrt{2}) = (\{x + y\sqrt{2} \ / \ x, y \in \mathbb{Q}\}, +, \cdot),$$

consisting of numbers of the form $x + y\sqrt{2}$ with the operation rules (1). In accordance with the dialectical law of negation of negation, we have again returned to what was found by the historical process, but at a higher level: now $\sqrt{2}$, $-\frac{1}{3}\sqrt{2}$, $\frac{3}{7} + 5\sqrt{2}$, are not some mystical "deaf numbers", but otherwise designated pairs $(0, 1), \left(0, -\frac{1}{3}\right), \left(\frac{3}{7}, 5\right)$ of the field $\widehat{\mathbb{Q}}$ or the corresponding elements of any field isomorphic to it.

In a similar way, one can obtain from the field of rational numbers an extension $\mathbb{Q}(\sqrt{3})$, in which the equation $x^2 = 3$ is solvable, also an extension $\mathbb{Q}(\sqrt[3]{2})$, in which $x^3 = 2$ is solvable, etc. Two successive extensions construct field $\mathbb{Q}(\sqrt{2}, \sqrt{3}) = \left(\mathbb{Q}(\sqrt{2})\right)(\sqrt{3})$, where both equations $x^2 = 2$ and $x^2 = 3$ are solvable, etc. By making reference to such examples in practical classes, we will shed light on the question of the "greatest economy" of extension, for which we will again use the field $\mathbb{Q}(\sqrt{2})$.

This is not the only extension of the field $\mathbb{Q}$ that leads to the solvability of the equation $x^2 = 2$: the same property is possessed, for example, by the field $\mathbb{Q}(\sqrt{2}, \sqrt{3})$, the field $\mathbb{R}$ of real numbers, etc. However, uniqueness (up to isomorphism) will take place if an additional *minimality* condition is imposed on the extension: in no proper subfield of this extension containing $\mathbb{Q}$ is the equation $x^2 = 2$ solvable.

**Theorem 11.** *All minimal extensions of the field of rational numbers that lead to the solvability of the equation $x^2 = 2$ are isomorphic to each other.*

*Proof.* The minimal extension of the field $\mathbb{Q}$ containing the root $\sqrt{2}$ of equation $x^2 = 2$ must also contain all numbers of the form $x + y\sqrt{2}$, where $x, y \in \mathbb{Q}$ (otherwise it will not even be a groupoid with respect to addition and multiplication). It is clear that if $x_1 = x_2$ and $y_1 = y_2$, then $x_1 + y_1\sqrt{2} = x_2 + y_2\sqrt{2}$; we will show that the opposite is also true.

Let $x_1 + y_1\sqrt{2} = x_2 + y_2\sqrt{2}$ be the same element of the extension. Then $(y_1 - y_2)\sqrt{2} = x_2 - x_1$. If $y_1 = y_2$, then the expression on the left side $(y_2 - y_1)\sqrt{2} = 0$ means that on the right side it must occur that $x_2 - x_1 = 0$ and so, $x_1 = x_2$. The statement is proven. The case $y_1 \neq y_2$ is impossible because then it would be $\sqrt{2} = \frac{x_2 - x_1}{y_1 - y_2} \in \mathbb{Q}$.

Due to the uniqueness of the representation of the numbers of the extended field in the form $x + y\sqrt{2}$ $(x, y \in \mathbb{Q})$, the mapping $f$ of this extension onto the field of pairs $\widehat{\mathbb{Q}}$ is defined as follows:
$$f(x + y\sqrt{2}) = (x, y),$$

which is a bijection. But bijection $f$ is an isomorphism:

$$f\left((x_1 + y_1\sqrt{2}) + (x_2 + y_2\sqrt{2})\right) = f\left((x_1 + x_2) + (y_1 + y_2)\sqrt{2}\right) = (x_1 + x_2, y_1 + y_2) =$$
$$= (x_1, y_1) + (x_2, y_2) = f(x_1 + y_1\sqrt{2}) + f(x_2 + y_2\sqrt{2}),$$

$$f\big((x_1 + y_1\sqrt{2})(x_2 + y_2\sqrt{2})\big) = f\big((x_1 x_2 + 2y_1 y_2) + (x_1 y_2 + y_1 x_2)\sqrt{2}\big) =$$
$$= (x_1 x_2 + 2y_1 y_2, x_1 y_2 + y_1 x_2) = (x_1, y_1)(x_2, y_2) = f(x_1 + y_1\sqrt{2}) \cdot f(x_2 + y_2\sqrt{2}).$$

And so, we have shown that any minimal extension of the field $\mathbb{Q}$ containing the root of equation $x^2 = 2$ is isomorphic to the field $\widehat{\mathbb{Q}}$ (or, equivalently, to the field $\mathbb{Q}(\sqrt{2})$). But since the relation of isomorphism of fields (as well as any algebraic systems in general) has the property of *symmetry* (reciprocity),

*if $P_1$ is isomorphic to $P_2$, then $P_2$ is isomorphic to $P_1$*

and the property of *transitivity*,

*if $P_1$ is isomorphic to $P_2$ and $P_2$ is isomorphic to $P_3$, then $P_1$ is isomorphic to $P_3$.*

We have already spoken about the first property of the isomorphism relation and we propose that you prove the second one, that is, any two fields that are isomorphic to the same third are isomorphic to each other. The theorem is proven. □

Using the field $\mathbb{Q}(\sqrt{2})$ as an example, we will examine another interesting question: which of the two roots of equation $x^2 = 2$ should we take as the number $\sqrt{2}$? The answer "positive root" is only valid when using such additional properties of numbers, which are expressed by the order relations $>, <$ in the field $\mathbb{Q}$ and, as we will see below, do not follow from the general definition of a field (and cannot be naturally introduced in every field). If we proceed only from conditions I, II and III, then we can only say about the two roots $\sqrt{2}$ and $-\sqrt{2}$ of the equation $x^2 = 2$ that they are opposite to each other, but calling the first of them "positive" and the second "negative" is meaningless: with equal success we could designate the pair $(0, -1)$ with the symbol $\sqrt{2}$, and then the pair $(0, 1)$ would receive the designation $-\sqrt{2}$.

A similar situation arises when trying somehow to distinguish the vertices at the base of an isosceles triangle while observing its specific location in the plane, which allows one of the vertices to be considered "left" and the other "right" (Fig.12). The complete equality (indistinguishability from each other) of both vertices is mathematically expressed in the fact that there is such a transformation (namely, a mirror image or reflection about its altitude), which carries the entire triangle onto itself, while the base vertices interchange (swap) positions. Similarly, the indistinguishability in the field $\mathbb{Q}(\sqrt{2})$ of the two roots of equation $x^2 = 2$ is expressed by the theorem proved below.

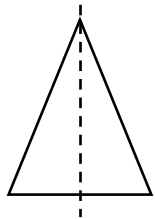

Fig. 12

An isomorphism of a field onto itself is called an *automorphism*. In addition to the *trivial automorphism* $f(x) = x$ that maps each element to itself, a field may admit other automorphisms.

**Theorem 12.** *Field* $\mathbb{Q}(\sqrt{2})$ *admits an automorphism* $f$ *such that* $f(\sqrt{2}) = -\sqrt{2}$ *(i.e., $f$ maps one root of equation* $x^2 = 2$ *to the other).*

*Proof.* We define a mapping $f: \mathbb{Q}(\sqrt{2}) \to \mathbb{Q}(\sqrt{2})$ of the field $\mathbb{Q}(\sqrt{2})$ onto itself assuming

$$f(x + y\sqrt{2}) = x - y\sqrt{2} \text{ for any } x, y \in \mathbb{Q}.$$

$f$ is a bijection (prove it!) and is also an automorphism because

$$f\big((x_1 + y_1\sqrt{2}) + (x_2 + y_2\sqrt{2})\big) = f\big((x_1 + x_2) + (y_1 + y_2)\sqrt{2}\big) =$$
$$= (x_1 + x_2) - (y_1 + y_2)\sqrt{2} = (x_1 - y_1\sqrt{2}) + (x_2 - y_2\sqrt{2}) = f(x_1 + y_1\sqrt{2}) + f(x_2 + y_2\sqrt{2}),$$
$$f\big((x_1 + y_1\sqrt{2})(x_2 + y_2\sqrt{2})\big) = f\big((x_1 x_2 + 2y_1 y_2) + (x_1 y_2 + y_1 x_2)\sqrt{2}\big) =$$
$$= (x_1 x_2 + 2y_1 y_2) - (x_1 y_2 + y_1 x_2)\sqrt{2} = (x_1 - y_1\sqrt{2})(x_2 - y_2\sqrt{2}) =$$
$$= f(x_1 + y_1\sqrt{2}) \cdot f(x_2 + y_2\sqrt{2}).$$

Automorphism $f$ is the desired one since in the case $x = 0, y = 1$, we obtain $f(\sqrt{2}) = -\sqrt{2}$. The theorem now follows. □

*Remark 15.* For $x = 0$, $y = -1$, we have $f(-\sqrt{2}) = \sqrt{2}$, that is, the automorphism $f$, which is discussed in Theorem 12, swaps the roots of the equation $x^2 = 2$. The example of the field $\mathbb{Q}(\sqrt{2})$ itself is perhaps not so important, but we used it as a convenient means of explaining the general idea of field extension.

Let $P$ be an arbitrary field, and equation

$$x^n + a_{n-1}x^{n-1} + \cdots + a_1 x + a_0 = 0,$$

where $n \in \mathbb{N}$ and all $a_{n-1}, \ldots, a_1, a_0 \in P$, has no roots in $P$. A *simple algebraic extension* $P(\theta)$ of field $P$ is its minimal extension in which the above equation of $n$th degree is solvable (by $\theta$, we mean any of the roots of this equation belonging to the extended field). It is curious that a non-simple algebraic extension obtained as a result of several simple ones can be represented as a single simple. We will not bring attention to more precise formulations and further study of this issue in the general course, but for those who are particularly interested, we will pose the following problem: find such an element $\theta$ and such an equation (with coefficients from $\mathbb{Q}$) that it satisfies, so that the simple algebraic extension $\mathbb{Q}(\theta)$ coincides with the field $\mathbb{Q}(\sqrt{2}, \sqrt{3})$ obtained from $\mathbb{Q}$ by two successive simple extensions.

At the same time, no algebraic extensions can transform $\mathbb{Q}$ into the field $\mathbb{R}$ of real numbers, if only because the latter has the *power of the continuum* (real numbers are uncountable), while $\mathbb{Q}$ itself and all fields obtained from it by successive simple extensions are only countable. Questions related to *cardinality* (a generalization of the concept of the number of elements applied to infinite sets) are initially considered in the course of mathematical analysis. The power set of a denumerable set is non-enumerable, and so its cardinality is larger than that of any denumerable set, which is $\aleph_0$ (aleph-null), standard for the cardinal of $\mathbb{N}$. The size of $\wp(\mathbb{N})$, where N is defined as the set of all subsets of N, is called the "power of the continuum" since it is the same size as the points on the real number line, $\mathbb{R}$. To obtain $\mathbb{R}$ from $\mathbb{Q}$, an analytical approach is used based on the axiom of completeness (or Dedekind cuts, fundamental sequences, etc.) and rooted in adding to $\mathbb{Q}$ such new numbers that allow their approximation by a rational number with any given accuracy. The latter makes sense not in any field, but, for example, in those whose elements can be compared with each other ("more than", "less than").

A field $(M, +, \cdot)$ is called an *ordered* one if a given binary order relation $\leq$ is defined on set $M$ satisfying the following conditions for any $x, y, z \in M$:

(1) $x \leq x$;
(2) if $x \leq y$ and $y \leq z$, then $x \leq z$;
(3) if $x \leq y$ and $y \leq x$, then $x = y$;
(4) always $x \leq y$ or $y \leq x$, that is, any two elements $x, y \in M$ are *comparable*;
(5) if $x \leq y$, then $x + z \leq y + z$;
(6) if $0 \leq x$ and $0 \leq y$, then $0 \leq xy$.

In particular, the field of rational numbers $\mathbb{Q}$ is an ordered one by the relation $\leq$ that has the generally accepted meaning. The extension (non-algebraic!) $\mathbb{R}$ of field $\mathbb{Q}$ is an ordered extension by the same relation.

Defining on an ordered field $(M, +, \cdot)$, which can now be written as $(M, +, \cdot, \leq)$, the relation $x < y$ as "$x \leq y$ and $x \neq y$", we can replace axioms (3) and (4) by one: for any $x, y \in M$ one and only one of the following three relations is true: $x = y$, $x < y$, $x > y$. After such a replacement, the ordered field itself can be designated by $(M, +, \cdot, <)$.

Not every field can be ordered. If, for example, we assume that it was possible to introduce order relations to the field $\mathbb{Z}_2$, then we will arrive at a contradiction: assuming first that $0 \leq 1$, based on (5), we obtain $1 = 0 + 1 \leq 1 + 1 = 0$ and further, due to (3), $0 = 1$ despite the fact that $\mathbb{Z}_2$ is a field. The assumption $1 \leq 0$ leads to a similar contradiction.

The field $\mathbb{R}$ of real numbers, unlike its "parent" $\mathbb{Q}$, has the property of completeness: in it, every Dedekind cut determines some number, or, equivalently, every fundamental sequence (satisfying the Cauchy condition) converges (has a limit). Another equivalent form is often used: any sequence of intervals in which each subsequent one belongs to the previous one has a common point (the *principle of nested intervals*). Algebraic extensions $\mathbb{Q}(\sqrt{2})$, $\mathbb{Q}(\sqrt{3})$, $\mathbb{Q}(\sqrt{2}, \sqrt{3})$ etc., are proper subfields of field $\mathbb{R}$; soon we will encounter such a simple algebraic extension of the field $\mathbb{Q}$ that goes beyond $\mathbb{R}$, and we will now indicate the following classification of real numbers.

A number $x \in \mathbb{R}$ is called *algebraic* if it satisfies some algebraic equation with rational coefficients, that is, if there exists such $n \in \mathbb{N}$ and such $a_0, a_1, \ldots, a_n \in \mathbb{Q}$, where $a_n \neq 0$, that

$$a_0 + a_1 x + a_2 x^2 + \cdots + a_n x^n = 0.$$

A number that is not algebraic, that is, does not satisfy any equation of the specified type, is called *transcendental*.

The above considerations on the cardinalities of sets show that there is a continuum of transcendental numbers, but they do not yet allow us to name even one of them specifically. Such examples were found in 1844 by J. Liouville, but these were numbers specially constructed by him, not previously known, and only in 1883 did Ch. Hermite prove the transcendence of the number $e = 2.7182818284\ldots$ (the base of natural logarithms), which by that time had already become firmly established in mathematics. The transcendence of the "truly ancient" number $\pi = 3.1415926535$ was established in 1892 by F. Lindemann, and in 1907, D. Hilbert greatly simplified both proofs. Based on Lindemann's theorem, one can, for example, immediately say without any calculations that the equality $\pi^3 - 10\pi + 0.41 = 0$ is not exact.

The set $\mathbb{Q}[\pi]$ of all numbers of the form $a_0 + a_1\pi + \cdots + a_n\pi^n, n = 0, 1, 2, \ldots;\ a_0, a_1, \ldots, a_n \in \mathbb{Q}$, forms in $\mathbb{R}$ only a subring, but not a subfield since already $\frac{1}{\pi} \notin \mathbb{Q}[\pi]$: if a representation $\frac{1}{\pi} = a_0 + a_1\pi + \cdots + a_n\pi^n$ had existed, then $-1 + a_0\pi + a_1\pi^2 + \cdots + a_n\pi^{n+1} = 0$ would have followed, contrary to Lindemann's theorem. But also adding to the sums $a_0 + a_1\pi + \cdots + a_n\pi^n$ a term of the form $\frac{b}{\pi}$, where $b \in \mathbb{Q}$, still does not transform $\mathbb{Q}[\pi]$ into a field: we suggest to show for yourself as an exercise that the number $\frac{1}{1+\pi}$ does not allow the representation $\frac{1}{1+\pi} = \frac{b}{\pi} + a_0 + a_1\pi + \cdots + a_n\pi^n$, where $b, a_0, a_1, \ldots, a_n \in \mathbb{Q}$. Looking ahead, we will remark that a *simple transcendental extension* $\mathbb{Q}(\pi)$, that is, a minimal field containing a subfield $\mathbb{Q}$ and the number $\pi$, consists of all possible numbers of the form $\frac{a_1 + a_2\pi^2 + \cdots + a_n\pi^n}{b_0 + b_1\pi + \cdots + b_m\pi^m}$, where $n, m \in \mathbb{N}_0$, $a_1, \ldots, a_n, b_0, \ldots, b_m \in \mathbb{Q}$ and $b_0, \ldots, b_m$ are not all zero.

## FIELD OF COMPLEX NUMBERS

As was already mentioned in the introductory lecture, the solution of 3rd and 4th degree equations in the 16th century led to the need to operate with numbers of an incomprehensible nature having the form $x + y\sqrt{-1}$, where the "imaginary" ("sophistical") number $\sqrt{-1}$ is such that $(\sqrt{-1})^2 = \sqrt{-1} \cdot \sqrt{-1} = -1$. Assuming that the three properties - commutative, associative and distributive - also remain valid for the new numbers, we obtain the rule of addition and multiplication:

$$\left.\begin{array}{l}(x_1 + y_1\sqrt{-1}) + (x_2 + y_2\sqrt{-1}) = (x_1 + x_2) + (y_1 + y_2)\sqrt{-1} \\ (x_1 + y_1\sqrt{-1}) \cdot (x_2 + y_2\sqrt{-1}) = (x_1 x_2 - y_1 y_2) + (x_1 y_2 + y_1 x_2)\sqrt{-1}\end{array}\right\},$$

where $x_1, x_2, y_1, y_2$ are any real numbers.

To substantiate this historical discovery, we, as in the case of adding $\sqrt{2}$ to the field of rational numbers, will construct such an extension of the field $\mathbb{R}$ of real numbers, in which equation $z^2 = -1$ is solvable. As before, the fact that the number $x + y\sqrt{-1}$ is completely determined by assigning an ordered pair $(x, y)$ of real numbers to it suggests the idea of taking all possible such pairs as elements of $\mathbb{R}$, and the above rules suggest the following definition of addition and multiplication of pairs:

$$(x_1, y_1) + (x_2, y_2) = (x_1 + x_2, y_1 + y_2),$$
$$(x_1, y_1) \cdot (x_2, y_2) = (x_1 x_2 - y_1 y_2, x_1 y_2 + y_1 x_2).$$

We will show that the system $(\{(x, y)/x, y \in \mathbb{R}\}, +, \cdot)$ is a field containing a subfield isomorphic to $\mathbb{R}$, and also that in this extended field the equation corresponding to $z^2 = -1$ is solvable. In this case, we will bring attention in detail only to those points that are essentially different from the points of reasoning when expanding the field $\mathbb{Q}$ to $\mathbb{Q}(\sqrt{2})$. The difference is due only to the fact that now $x$ and $y$ are any real numbers (not necessarily rational), and in the definition of multiplication of pairs; the first component of the product has the form $x_1 x_2 - y_1 y_2$ (instead of $x_1 x_2 + 2 y_1 y_2$). Let's introduce the notation $\widehat{\mathbb{R}} = \{(x, y)/x, y \in \mathbb{R}\}$ and recall once again that the equality of pairs $(x_1, y_1) = (x_2, y_2)$ means $x_1 = x_2$ and $y_1 = y_2$.

The fact that $(\widehat{\mathbb{R}}, +)$ and $(\widehat{\mathbb{R}}, \cdot)$ are groupoids, moreover, commutative and associative, and that in $(\widehat{\mathbb{R}}, +, \cdot)$, addition and multiplication satisfy the distributive property, can be proved according to the previous scheme and the reader should be able to do this independently. We just need to refer to the same laws in $\mathbb{Q}$ for commutativity, associativity and distributivity in $\mathbb{R}$, and also not forget about the difference in the rule for multiplying pairs. The conclusions that $(\widehat{\mathbb{R}}, +)$ is an abelian group and $(\widehat{\mathbb{R}}, +, \cdot)$ is an associative-commutative ring with unity $(1, 0)$ also remain valid; let's check the latter: $(x, y) \cdot (1, 0) = (1, 0) \cdot (x, y) = (1 \cdot x - 0 \cdot y, 1 \cdot y + 0 \cdot x) = (x, y)$.

To prove the invertibility of any element $(a, b) \neq (0, 0)$ from $\widehat{\mathbb{R}}$, we first appeal to the "historical" transformation

$$\frac{1}{a+b\sqrt{-1}} = \frac{a-b\sqrt{-1}}{(a+b\sqrt{-1})(a-b\sqrt{-1})} = \frac{a-b\sqrt{-1}}{a^2 - b^2(\sqrt{-1})^2} = \frac{a-b\sqrt{-1}}{a^2 + b^2} = \frac{a}{a^2 + b^2} - \frac{b}{a^2 + b^2}\sqrt{-1},$$

where $a^2 + b^2 \neq 0$ (the sum of the squares of two real numbers, of which at least one is different from zero). This suggests the idea of immediately checking that the pair $\left(\frac{a}{a^2 + b^2}, \frac{-b}{a^2 + b^2}\right)$ serves as an inverse element for $(a, b)$. And indeed,

$$(a, b) \cdot \left(\frac{a}{a^2 + b^2}, \frac{-b}{a^2 + b^2}\right) = \left(a \cdot \frac{a}{a^2 + b^2} - b \cdot \frac{-b}{a^2 + b^2}, a \cdot \frac{-b}{a^2 + b^2} + b \cdot \frac{a}{a^2 + b^2}\right) =$$
$$= \left(\frac{a^2 + b^2}{a^2 + b^2}, \frac{-ab + ba}{a^2 + b^2}\right) = (1, 0).$$

In view of the commutativity of the multiplication of pairs, it is unnecessary to perform a check with permuted factors (although it is not bad as an exercise), and the uniqueness of the inverse element was established earlier in the general case (see Theorem 2).

The end of the proof that $(\widehat{\mathbb{R}}, +, \cdot)$ is a field remains literally the same. There is no difference in the proof that the pairs of the form $(x, 0)$ form a subfield isomorphic to $\mathbb{R}$ with a natural isomorphism $f(x) = (x, 0)$. In particular, $f(-1) = (-1, 0)$. We will not only prove the solvability in $\widehat{\mathbb{R}}$ of the equation $(x, y)^2 = (-1, 0)$, but also find all its solutions: now this question is more important than in the construction of $\mathbb{Q}(\sqrt{2})$.

The equation
$$(x, y)^2 = (x, y) \cdot (x, y) = (-1, 0)$$
in expanded form
$$(x^2 - y^2, 2xy) = (-1, 0)$$
is equivalent to the system of two equations
$$x^2 - y^2 = -1,$$
$$2xy = 0,$$
for which all real solutions must be found.

Due to the absence of zero divisors in $\mathbb{R}$, the second equation means that $x = 0$ or $y = 0$. For $x = 0$, from the first equation, we get $y^2 = 1$, that is, $y = \pm 1$, and for $y = 0$, the first equation takes the form $x^2 = -1$ that has no real solutions. Therefore, the pairs $(0, 1), (0, -1)$ and only these two pairs are solutions of the equation $(x, y)^2 = (-1, 0)$.

Considering the symbol of the number $x \in \mathbb{R}$ as another designation of the pair $(x, 0) \in \widehat{\mathbb{R}}$ and denoting the pair $(0, 1)$ by the letter $i$, we obtain for any element $(x, y) \in \widehat{\mathbb{R}}$ the representation

$$(x, y) = (x, 0) + (y, 0)(0, 1) = x + y \cdot i = x + iy$$

(it is customary to write $x + iy$, but at the same time $a + bi$ rather than $a + ib$).

So, the desired extension $\mathbb{R}(\sqrt{-1})$ of the field $\mathbb{R}$ consists of numbers of the form $x + iy$, where $x, y \in \mathbb{R}$ and $i^2 = -1$, and these are no longer "mystical" numbers, but only differently designated pairs of real numbers (remember again the dialectical law of negation of negation, similar to the case of $\mathbb{Q}(\sqrt{2})$, is true).

**Theorem 13.** *All minimal extensions of the field of real numbers that lead to the solvability of the equation $z^2 = -1$ are isomorphic to each other.*

*Proof.* The minimal extension of the field $\mathbb{R}$ containing the root $i$ of equation $z^2 = -1$ must also contain all numbers of the form $x + iy$ with $x, y \in \mathbb{R}$. We will show that if $x_1 + iy_1 = x_2 + iy_2$, then $x_1 = x_2$ and $y_1 = y_2$.
As with the proof of Theorem 11, let $x_1 + iy_1 = x_2 + iy_2$ be the same element of the extension. Then $i(y_1 - y_2) = x_2 - x_1$ and either $y_1 = y_2 \Rightarrow x_1 = x_2$, (as required) or $y_1 \neq y_2$. But the second case is impossible because then it would be $i = \frac{x_2 - x_1}{y_1 - y_2} \in \mathbb{R}$.

Due to the uniqueness of the representation of the numbers of the extended field in the form of $x + iy$, the mapping
$$f(x + iy) = (x, y)$$

is a bijection of the minimal extension $\mathbb{R}(i)$ under consideration onto the field $\mathbb{R}$ of pairs $\widehat{\mathbb{R}}$. This bijection is an isomorphism:

$$f\big((x_1 + iy_1) + (x_2 + iy_2)\big) = f\big((x_1 + x_2) + i(y_1 + y_2)\big) = (x_1 + x_2, y_1 + y_2) =$$
$$= (x_1, y_1) + (x_2, y_2) = f(x_1 + iy_1) + f(x_2 + iy_2),$$
$$f\big((x_1 + iy_1)(x_2 + iy_2)\big) = f\big((x_1 x_2 - y_1 y_2) + i(x_1 y_2 + y_1 x_2)\big) =$$
$$= (x_1 x_2 - y_1 y_2, x_1 y_2 + y_1 x_2) = (x_1, y_1)(x_2, y_2) = f(x_1 + iy_1) \cdot f(x_2 + iy_2).$$

In a similar way, as was done with the proof of Theorem 11, it follows that all minimal extensions of the field $\mathbb{R}$ leading to the solvability of equation $z^2 = -1$ are isomorphic to field $\widehat{\mathbb{R}}$, and therefore to each other. The statement now follows. □

A minimal extension $\mathbb{R}(i)$ of field $\mathbb{R}$ containing element $i$, for which $i^2 = -1$, is called the *field of complex numbers* and is denoted by $\mathbb{C}$. Any field isomorphic to it can also be called and denoted in the same way, for example $\widehat{\mathbb{R}}$.

A non-trivial automorphism $f(x + iy) = x - iy$ of field $\mathbb{C}$ is called the *conjugate* of $f$ and it is designated using a bar over the letter $f$: if $z = x + iy \in \mathbb{C}$, then its complex conjugate, by definition, is

$$\overline{z} = \overline{x + iy} = x - iy.$$

The properties of conjugates of complex numbers, the geometric representation of complex numbers, and also, their trigonometric form will be presented traditionally (see A.G. Kurosh, Course of Higher Algebra). *De Moivre's formula*

$$z^n = |z|^n (\cos n\varphi + i \sin n\varphi),$$

where $\varphi = \arg(z) =$ the angle between the positive real axis and the line joining the origin and $z$, $|z| = \sqrt{x^2 + y^2}$, is proven by induction on $n \in \mathbb{N}$.

The presentation of the topic "extracting the $n$th root of a complex number" differs from the traditional one only in that when assigning a set

$$\{z \in \mathbb{C} \,/\, z^n = c\} = \sqrt[n]{|c|}\left[\cos \frac{\arg(c) + 2\pi k}{n} + i \sin \frac{\arg(c) + 2\pi k}{n}\right] \qquad (2)$$

$k = 0, 1, \ldots, n - 1$, we do not denote it by $\sqrt[n]{c}$ for the following reason: using the notation

$$\sqrt[n]{c} = \sqrt[n]{|c|}\left[\cos \frac{\arg(c) + 2\pi k}{n} + i \sin \frac{\arg(c) + 2\pi k}{n}\right],$$

the symbol $\sqrt[n]{\phantom{c}}$ on the left and on the right would have different meanings. On the right; an ordinary arithmetic root (a positive number with $c \neq 0$), and on the left; a "multi-valued function", which is why with a real $c > 0$ we get an absurdity, for example

$$\sqrt[3]{8} = \sqrt[3]{8}\left(\cos\frac{2\pi k}{3} + i\sin\frac{2\pi k}{3}\right),$$

and reducing both sides by $\sqrt[3]{8}$, we arrive at an equality that is not true for $k = 1, 2$. For "multi-valued functions" even the simplest relations like $\sqrt[3]{8} + \sqrt[3]{8} = 2\sqrt[3]{8}$ are not true (the left side has six different values and the right side has only three as shown in Fig.13).

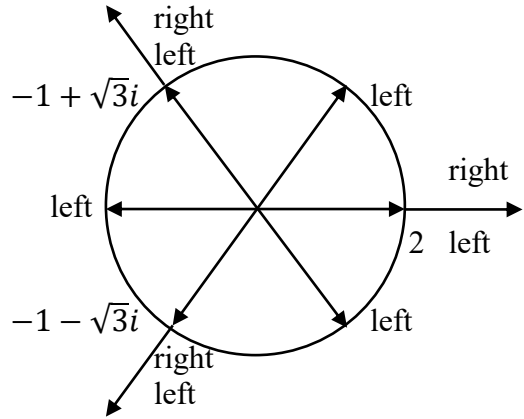

Fig. 13

Paying tribute to the deep-rooted tradition, we will sometimes use the sign $\sqrt[n]{c}$ for a complex $c$, but each time with limitations. Note that even in problem books there still appears an imprecise transcription like "find all the values of $\sqrt[4]{i}$".

*The field $\mathbb{C}$ of complex numbers cannot be ordered*, that is, a binary relation $\leq$ cannot be defined on the pairs of its elements that would satisfy axioms (1) – (6) of an ordered field. Indeed, assume this is possible. Since $i \neq 0$, by axiom (4), either $i < 0$ or $0 < i$. If $i < 0$, then by axiom (5), $0 \leq -i$ and by axiom (6), $0 \leq (-i)(-i)$, that is, $0 \leq -1$; from here and from $0 \leq -i$, we get, again by axiom (6), $0 \leq i$, therefore $0 < i$ (since $i \neq 0$), contrary to the assumption. A similar contradiction arises assuming $0 < i$.

SOLVING ALGEBRAIC EQUATIONS

<u>1st degree equation</u> $\alpha z + \beta = 0$, where $\alpha, \beta \in \mathbb{C}$ and $\alpha \neq 0$ (case $\alpha = 0$ is not even an equation, but rather the banal identity $0 = 0$ or a false statement like $1 = 0$). Obviously, the solution is unique: $z = -\frac{\beta}{\alpha}$. But even such simple things required from us a new justification since we now are dealing with the field of complex numbers, to which, mechanically, it is not possible to transfer what was previously known only for real numbers.

<u>2nd degree equation</u> $\alpha z^2 + \beta z + \gamma = 0$, where $\alpha, \beta, \gamma \in \mathbb{C}, \alpha \neq 0$ (otherwise, the equation is not of second degree). Dividing both sides by $\alpha$ and introducing the notation $p = \frac{\beta}{\alpha}, q = \frac{\gamma}{\alpha}$, we obtain the *reduced quadratic equation*

$$z^2 + pz + q = 0,$$

equivalent to the original, where $p, q \in \mathbb{C}$.

To solve, we complete the square (of course, adding the same amount to the right) and perform other obvious transformations:

$$z^2 + 2z \cdot \frac{p}{2} + \left(\frac{p}{2}\right)^2 + q = \left(\frac{p}{2}\right)^2 \Rightarrow \left(z + \frac{p}{2}\right)^2 = \frac{p^2}{4} - q,$$

and if by $\sqrt{\frac{p^2}{4} - q}$, we mean the value of the root that corresponds to the value $k = 0$ (see formula (2)) for a set $\{w \;/\; w^2 = \frac{p^2}{4} - q\}$ (of two elements), then formula

$$z = -\frac{p}{2} \pm \sqrt{\frac{p^2}{4} - q}$$

will be correct since the $n$th root with $n = 2$ has a value for $k = 1$ that is opposite to the value for $k = 0$.

The *reduced quadratic trinomial* on the left side of the reduced quadratic equation is obviously a perfect square trinomial if and only if $q = \frac{p^2}{4}$. Expressing here $p$ and $q$ through $\alpha, \beta, \gamma$, after some simple transformations, we obtain a similar condition in the general case:

*a square trinomial $\alpha z^2 + \beta z + \gamma$ is a perfect square trinomial, that is, it has the form $(\delta z + \varepsilon)^2$ with $\delta, \varepsilon \in \mathbb{C}$ if and only if $\beta^2 = 4\alpha\gamma$.*

<u>3rd degree equation</u> $\alpha z^3 + \beta z^2 + \gamma z + \delta = 0;\ \alpha, \beta, \gamma, \delta \in \mathbb{C};\ \alpha \neq 0$. Making the change of variable $z = \tilde{z} - \frac{\beta}{3\alpha}$ and dividing both parts by $\alpha$, after substitution, some transformations and the already familiar designations of the coefficients, and also, omitting the sign $\sim$ over the new unknown, we obtain the *reduced cubic equation*

$$z^3 + pz + q = 0,$$

not containing the square of the unknown. The method for solving it was discovered by del Ferro (for a number of special cases), Tartaglia (in general form), and Cardano wrote down the general formula and researched all possibilities.

We will look for a solution in the form $z = u + v$, so that we can then get rid of the prevailing alternative of choosing $u$ and $v$.

$$(u + v)^3 + p(u + v) + q = 0,$$
$$u^3 + v^3 + (3uv + p)(u + v) + q = 0.$$

We select $u$ and $v$ so that $3uv + p = 0$. As a result, for the two unknowns $u^3$ and $v^3$, we obtain a system of two equations
$$u^3 + v^3 = -q$$
$$u^3 \cdot v^3 = -\frac{p^3}{27},$$
the second of which is the result of raising to the power of three both sides of equation $uv = -\frac{p}{3}$.

Knowing the sum and the product of two quantities ($u^3$ and $v^3$), which is the case when we are factoring a simple square trinomial, we arrive at the conclusion that they are the solution of the quadratic equation in factored form
$$(w - u^3)(w - v^3) = 0.$$
Expanding and grouping
$$w^2 - (u^3 + v^3)w + u^3 \cdot v^3 = 0,$$
we obtain the following *auxiliary quadratic equation*:
$$w^2 + qw - \frac{p^3}{27} = 0,$$
from which
$$w = -\frac{q}{2} \pm \sqrt{\frac{q^2}{4} + \frac{p^3}{27}},$$
that is,
$$u^3 = -\frac{q}{2} + \sqrt{\frac{q^2}{4} + \frac{p^3}{27}}, \quad v^3 = -\frac{q}{2} - \sqrt{\frac{q^2}{4} + \frac{p^3}{27}}$$

(or vice versa). From here, taking the cubic root of both sides, we obtain for the unknown $z = u + v$, the expression in the traditional imprecise notation:

$$z = \sqrt[3]{-\frac{q}{2} + \sqrt{\frac{q^2}{4} + \frac{p^3}{27}}} + \sqrt[3]{-\frac{q}{2} - \sqrt{\frac{q^2}{4} + \frac{p^3}{27}}}$$

(called Cardano's formula), where for a pair of terms $u, v$, all $3 \cdot 3 = 9$ values are taken and only three such that $uv = -\frac{p}{3}$.

Omitting a detailed analysis of all possible cases, we will only remark that if $p, q \in \mathbb{R}$ and the equation has three different real roots, then in Cardano's formula both terms turn out to be complex and the imaginary parts are cancelled out only as a result of addition. It also turns out to be impossible to obtain a similar formula that does not lead to the advent of imaginary numbers at the intermediate stages of the solution.

<u>4th degree equation</u> $\alpha z^4 + \beta z^3 + \gamma z^2 + \delta z + \varepsilon = 0$; $\alpha, \beta, \gamma, \delta, \varepsilon \in \mathbb{C}$; $\alpha \neq 0$. Dividing by $\alpha$, assuming $z = \tilde{z} - \frac{\beta}{4\alpha}$ and omitting the sign $\sim$, we obtain after some transformations the *reduced quartic equation*

$$z^4 + pz^2 + qz + r = 0,$$

which was solved by Ferrari, a student of Cardano.

Introducing an additional auxiliary unknown $w$, we write the equation in the form

$$\left(z^2 + \frac{p}{2} + w\right)^2 = 2w \cdot z^2 - q \cdot z + \left(w^2 + pw + \frac{p^2}{4} - r\right)$$

and we try to choose the value of $w$ so that the square trinomial (relative to $z$), standing on the right side, turns out to be a perfect square trinomial. For this there must be

$$(-q)^2 = 4 \cdot 2w \cdot \left(w^2 + pw + \frac{p^2}{4} - r\right);$$

this is an *auxiliary cubic equation* with respect to $w$, and we already have solutions for such equation. It is enough to find one of its roots $w_0$. Then, after the appropriate designation of the coefficients, the original equation will take the form

$$\left(z^2 + \frac{p}{2} + w_0\right)^2 = (\zeta z + \eta)^2,$$

wherefrom

$$z^2 + \frac{p}{2} + w_0 = \pm(\zeta z + \eta),$$

that is, the solution of the quartic equation was reduced to the solution of two quadratic ones.

"Taking the square root of both sides" is best justified as follows: Let $u^2 = v^2$. Then

$$u^2 - v^2 = 0 \Rightarrow (u - v)(u + v) = 0$$

wherefrom, due to the absence of zero divisors in $\mathbb{C}$, $(u - v) = 0$ or $(u + v) = 0$, that is,

$$u = \pm v.$$

Thus, in the 17th century, the algebraic equations of degree $n = 1, 2, 3, 4$ (and also the *binomial equation* $z^n = c$ of any degree) were solved. By solution was meant a *solution in radicals*, that is, a general formula (or a general method) for expressing the roots of an algebraic equation through its coefficients using the following operations:

(1) addition, subtraction, multiplication and division;
(2) extraction of the $n$th root, for any $n \in \mathbb{N}, n \geq 2$;
(3) selection of values from a given finite set.

In the 18th century, there was a search for solutions (in the same sense) of algebraic equations of any degree. An auxiliary equation, whose solution allows to solve the original one, is called its *resolvent*. Since the resolvent of the 4th degree equation itself has degree 3 and the resolvent of the 3rd degree equation has degree 2 (when solving a 2nd degree equation, it is not difficult to see a resolvent of degree 1), one could expect that for the 5th degree equation, the resolvent would be of degree 4, and so on. But the harsh reality turned out to be different: the resolvent of the 5th degree equation, found by Lagrange, has degree 6 (out of the frying pan into the fire!).

A proposal arose that the general algebraic equation of degree $n \geq 5$ is unsolvable in radicals. Ruffini gave an incomplete proof of this statement. The proof was completed (already at the beginning of the 19th century) by Abel, who is considered the founder of group theory (since he considered only commutative groups, the latter are now called abelian).

But the absence of a general solution to, let's say, a 5th degree equation in radicals does not mean that it is impossible to write such a solution for a specific equation (that is, with specific numerical values of the coefficients). For example, if the left-hand side can be factored. Developing group theory further, on its basis, Galois established a criterion for the solvability of specific equations in radicals (for instance, equation $z^5 - z - 1 = 0$ is unsolvable). *Galois theory* as a separate special course is taught in some universities and pedagogical institutes. We recommend reading the following popular books:
O. Ore, Niels Henrik Abel: Mathematician Extraordinary.
L. Infeld, Whom the Gods Love: The Story of Evariste Galois.

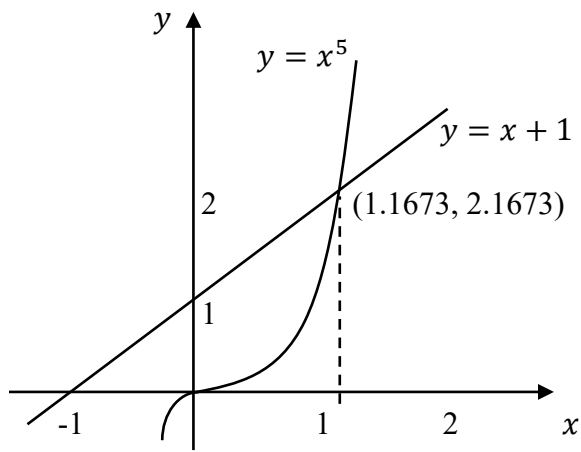

Fig. 14

In turn, the insolvability of an equation in radicals does not mean that it has no roots. Thus, the equation $x^5 = x + 1$ (Fig. 14) obviously has a real root (a little more than 1, which can be found approximately with any required accuracy: $x = 1.1673$ works well), this only root cannot be written as a result of actions (1), (2) and (3) over the coefficients $1, -1, -1$ of the equation. Moreover, the following holds:

**Theorem 14.** FUNDAMENTAL THEOREM OF ALGEBRA. *Every algebraic equation with coefficients belonging to the field $\mathbb{C}$ has at least one root in this field.* □

*Remark 16.* The proof of this theorem, begun by d'Alembert, was completed by Gauss. Later (19th century) other proofs appeared, but there is no a purely algebraic one among them and, as has now been established, there cannot be one (however, this is not surprising: after all, the field $\mathbb{C}$ is an extension of the field $\mathbb{R}$ of real numbers, constructed from $\mathbb{Q}$ not algebraically, but analytically). For this reason, we accept the fundamental theorem of algebra in the algebra course itself without a proof; it is easily obtained in the theory of functions of a complex variable as a consequence of Liouville's theorem (which we will not even formulate here and will only remark that its proof does not rely either directly or indirectly on the fundamental theorem of algebra, otherwise it would be a vicious circle!)

The fundamental theorem of algebra states that a further algebraic extension of the field $\mathbb{C}$ is impossible. It turns out that other extensions (non-algebraic or *transcendental*) are also impossible.

As we know, the extension of the field $\mathbb{R}$ to $\mathbb{C}$ is associated with the loss of order. By sacrificing the commutativity of multiplication, Hamilton constructed the skew field of quaternions containing $\mathbb{C}$ and consisting of numbers of the form $a + bi + cj + dk$ where $a, b, c, d \in \mathbb{R}$, and the "imaginary units" $i, j, k$ are multiplied by each other according to the law

$$i \cdot j = k, \quad j \cdot k = i, \quad k \cdot i = j, \quad j \cdot i = -k, \quad k \cdot j = -i, \quad i \cdot k = -j$$

As Frobenius showed, no other skew field (sometimes called *sfield*) containing $\mathbb{C}$ as a subfield exists (further "extensions" are associated with the rejection of such properties as, for example, associativity of multiplication, distributive laws, etc.).